%% file: main_elsevier_latex.tex
\documentclass[final,3p,times]{elsarticle}

\usepackage{iftex}

\ifPDFTeX
  \usepackage[utf8]{inputenc}
  \usepackage[T1]{fontenc}
  \usepackage{CJKutf8}
\else
  \ifLuaTeX
    \usepackage{luatexja}
    \usepackage{luatexja-fontspec}

    \setmainjfont{HaranoAjiMincho-Regular.otf}
    \setsansjfont{HaranoAjiGothic-Regular.otf}
  \else
    \errmessage{Please compile this file with LuaLaTeX or pdfLaTeX.}
  \fi
\fi

\usepackage{amsmath,amssymb}
\usepackage{graphicx}
\usepackage{xcolor}
\usepackage{multirow,array}
\usepackage{placeins}
\usepackage{comment}
\usepackage{subcaption}
\usepackage{algorithm}
\usepackage{algpseudocode}
\usepackage{makecell}
\usepackage[hidelinks]{hyperref}

\biboptions{sort&compress}

\newcommand{\bm}{\boldsymbol}

\DeclareFontShape{OMX}{cmex}{m}{n}{%
  <-> sfixed * [9] cmex10
}{}

\newcommand{\BeginJapaneseForPdfLaTeX}{%
  \ifPDFTeX
    \begin{CJK}{UTF8}{min}%
  \fi
}

\newcommand{\EndJapaneseForPdfLaTeX}{%
  \ifPDFTeX
    \end{CJK}%
  \fi
}

\journal{Elsevier journal name}

\begin{document}
\BeginJapaneseForPdfLaTeX

\begin{frontmatter}

\title{Metagraph-Based Domain-Decomposed\\ Galerkin Reduced-Order Model}

\author[aff1]{Kyohei Shintate}
\author[aff2]{Naoki Morita}
\author[aff3]{Shigeki Kaneko}
\author[aff1]{Nozomi Magome}
\author[aff2]{Naoto Mitsume\corref{cor1}}

\address[aff1]{University of Tsukuba, Degree Programs in Systems and Information Engineering, Tennodai 1-1-1, Tsukuba, 3058573, Ibaraki, Japan}
\address[aff2]{University of Tsukuba, Institute of Systems and Information Engineering, Tennodai 1-1-1, Tsukuba, 3058573, Ibaraki, Japan}
\address[aff3]{Graduate School of Engineering, Nagoya Institute of Technology, Gokishocho, Showa-ku, Nagoya, 4660061, Aichi, Japan}

\cortext[cor1]{Corresponding author. Email: mitsume@kz.tsukuba.ac.jp; Telephone: +81-29-853-5268}

\input{include/abstract}

\end{frontmatter}

\input{include/introductions}
\input{include/equations}
\input{include/numerical_analysis}
\input{include/conclusions}

\section*{CRediT authorship contribution statement}
\textbf{Kyohei Shintate}: Methodology, Software, Validation, Formal analysis, Investigation, Data curation, Writing -- original draft, Visualization, Funding acquisition.
\textbf{Naoki Morita}: Methodology, Software, Writing -- review \& editing, Supervision, Funding acquisition.
\textbf{Shigeki Kaneko}: Methodology, Writing -- review \& editing, Supervision.
\textbf{Nozomi Magome}: Software, Writing -- review \& editing, Funding acquisition.
\textbf{Naoto Mitsume}: Conceptualization, Project administration, Supervision, Writing -- review \& editing, Funding acquisition.

\section*{Funding}
This work was supported by JST FOREST Grant Number JPMJFR215S, JST ACT-X Grant Number JPMJAX24LG, Joint Usage/Research Center for Interdisciplinary Large-scale Information Infrastructures (JHPCN) Grant Number jh240017, and JSPS KAKENHI Grant Numbers 24H00695, 26K02921, 26K21236, 26KJ0663.

\section*{Declaration of conflict of interest}
The authors declare that they have no known competing financial interests or personal
relationships that could have appeared to influence the work reported in this paper.

\input{include/appendixA}
\input{include/appendixB}

\bibliographystyle{unsrt}
\bibliography{bibliography.bib}

\EndJapaneseForPdfLaTeX
\end{document}

%% file: include/abstract.tex
\renewcommand{\abstractname}{Abstract}

\begin{abstract}
This study proposes a metagraph-based domain-decomposed Galerkin reduced-order model (MBDD-G-ROM) for distributed-memory parallel reduced-order analysis of large-scale problems. The method provides a graph-based representation of domain-decomposed Galerkin reduced-order models defined over arbitrary domain decompositions. Calculation-point graphs describe interactions among calculation points, whereas metagraphs describe connectivity among local approximation-space subdomains. In the proper orthogonal decomposition (POD)-based implementation considered in this study, these local subdomains correspond to POD computation subdomains, where local POD bases are constructed. In the resulting metagraph, POD computation subdomains are treated as metanodes, and metaedges encode the block-sparsity pattern induced by overlaps between the supports of local POD basis functions. Based on this representation, a two-level graph-based domain-partitioning strategy decouples the POD computation subdomains from the parallel computation subdomains. This decoupling enables distributed-memory parallelization of the offline and online phases, including reduced-system assembly and the linear solver, without imposing a one-to-one correspondence between the two subdomain types in number or geometry. The metagraph structure also provides a natural way to incorporate static load balancing through metanode weights representing estimated computational costs. The method is verified using an unsteady diffusion equation and the incompressible Navier--Stokes equations for flow around a three-dimensional cylinder; these tests assess accuracy degradation caused by model reduction and strong-scaling behavior. Numerical results show that the method maintains solution accuracy while achieving high parallel efficiency in the online phase. The static load-balancing test shows that metanode weights based on estimated computational costs improve computational efficiency in the tested case.
\end{abstract}

\begin{keyword}

Parallel Computing \sep Domain Decomposition \sep Graph Structure \sep Reduced-Order Modeling \sep Proper Orthogonal Decomposition

\end{keyword}

%% file: include/introductions.tex
\section{Introduction}\label{sec:introduction}

Parametric studies play an important role in many engineering applications, where analyses are repeatedly performed for different values of the parameters defining the governing partial differential equation model. However, repeated simulations using full-order models (FOMs) entail substantial computational costs and are therefore often impractical. At the same time, because the target systems in such studies often share similar geometries and physical conditions across parameter instances, previously obtained simulation data can potentially be exploited to improve computational efficiency.
Recently, reduced-order models (ROMs) \citep{rom_book_vol2,rom_intro_book,LU2019264} have been studied extensively as fast surrogate models for FOMs. Among them, Galerkin reduced-order models (G-ROMs), which compute approximate solutions by projecting the governing equations onto a low-dimensional approximation space, either a linear subspace or a nonlinear manifold, have attracted considerable attention owing to their favorable balance between computational accuracy and efficiency.
Approaches to defining the low-dimensional approximation spaces used in G-ROMs can be broadly divided into two classes.
The first class employs prescribed linear subspaces, as exemplified by spectral methods \citep{CANUTO2006}, including spaces spanned by eigenmodes \citep{BENNIGHOF20042084} or trigonometric basis functions \citep{MOEHLIS200456,WALEFFE1997883}.
In this class of methods, the linear subspace is often chosen as the span of analytically defined functions, and such methods have been applied mainly to steady problems and problems with relatively simple geometries.
The second class consists of data-driven methods that construct a low-dimensional approximation space from FOM snapshots so as to represent the snapshot data as accurately as possible.
These data-driven approaches can be further categorized according to the form of the resulting approximation space into (A) linear-subspace-based approaches \citep{SIROVICH1987,WILLCOX20022323} and (B) nonlinear-manifold-based approaches \citep{RUTZMOSER2017196,BARNETT2022111348,GEELEN2023115717,KIM2024116978,LEE2020108973}.
In linear-subspace-based approaches, data-driven dimensionality-reduction techniques are applied to snapshot data to identify a basis spanning a linear subspace, and the solution field is approximated as a linear combination of the resulting basis vectors. This framework can achieve good accuracy even for problems involving complex geometries and solution fields. Representative methods include snapshot proper orthogonal decomposition (snapshot POD) \citep{SIROVICH1987} and balanced POD \citep{WILLCOX20022323}.
Representative examples of nonlinear-manifold-based approaches include machine-learning models such as autoencoders \citep{KIM2024116978,LEE2020108973}. These approaches are well suited to capturing strong nonlinearities and localized structures that are difficult to represent using linear subspaces. However, for problems with relatively weak nonlinearity, they can incur higher computational costs than linear-subspace-based approaches.
Such data-driven G-ROMs are generally organized into two stages: an offline phase and an online phase. In the offline phase, FOM simulations are performed for selected parameter samples, and a low-dimensional approximation space is constructed from the resulting snapshots. In the online phase, the governing equations are projected onto the identified low-dimensional approximation space, typically by Galerkin projection, thereby substantially reducing the number of unknowns and enabling fast computation. For example, the POD-Galerkin method \citep{LU2019264}, which combines snapshot POD \citep{SIROVICH1987} with Galerkin projection, has been applied to a wide range of engineering problems, including shape optimization \citep{LASSILA20101583,MANZONI2012}, control \citep{Ravindran2000,Strazzullo2020}, real-time analysis \citep{SON2013818,AMSALLEM20091241}, fluid dynamics \citep{BAIGES201323,CAO2019679,WANG2016374,EIXIMENO2025109459,FARCAS2025109619,CHUNG2024117041,peherstorfer2020model,NOACK2014,LUCIA2003917_snapshotPOD,BALLARIN2015,REYES2020112844,DAVE2025118194}, and structural mechanics \citep{KERFRIDEN2013169,CORIGLIANO2015127,ROCHA2020}.

As discussed above, data-driven ROMs can substantially accelerate simulations and reduce the number of unknowns.
However, memory requirements remain a major challenge in both the offline and online phases when the underlying FOM discretization involves an extremely large number of degrees of freedom (DOFs), as in direct numerical simulations of high-Reynolds-number flows.
For such complex problems, constructing a G-ROM requires training data generated by large-scale FOM simulations.
These training data consist of full-order solution vectors.
Moreover, the basis or mapping that defines the low-dimensional approximation space also contains quantities whose size scales with the number of full-order DOFs.
Consequently, storing, communicating, and manipulating the data needed to connect the full-order and reduced representations can become a critical bottleneck.
An effective way to address these issues is to use parallel computing on distributed-memory systems.
A representative approach is domain decomposition \citep{DOLEAN2015,DDMbook,FARHAT2001,DOHRMANN2003}, in which the cells or finite elements of the target domain are geometrically partitioned to distribute the data required for the analysis and balance the computational load.
Domain decomposition methods are widely used in numerical analysis and computational mechanics, particularly in finite element methods.
For POD-based G-ROMs, several studies have reported distributed-memory parallel computations using domain decomposition for the singular value decomposition (SVD) in the offline phase \citep{CAO2019679,WANG2016374,EIXIMENO2025109459,FARCAS2025109619}.
However, to the best of the authors' knowledge, a POD-based G-ROM framework with distributed-memory parallelization of the entire online phase, from the assembly of the reduced system of linear equations to its solution by a linear solver, has not yet been reported.
If the linear solver is not parallelized, the cost of solving the reduced linear system becomes non-negligible as the reduced dimension increases, leading to a loss of computational efficiency.
This issue is particularly relevant for advection-dominated problems \citep{NOACK2014} and turbulent-flow problems \citep{peherstorfer2020model}, which tend to require a large reduced dimension \citep{GALLETTI2004161}.
Furthermore, in G-ROMs that use a global low-dimensional approximation space constructed over the entire computational domain, the reduced coefficient matrix is generally dense.
Consequently, even if domain decomposition is applied to parallelize the reduced system of linear equations, interprocess communication tends to dominate the computation, which can lead to degraded parallel efficiency.

A different use of domain decomposition has also been proposed in the form of domain-decomposed Galerkin reduced-order models (DD-G-ROMs), where the decomposition is used not primarily for distributed-memory parallel computing but to construct local low-dimensional approximation spaces in individual subdomains.
In DD-G-ROMs, local low-dimensional approximation spaces are constructed independently from FOM data restricted to each subdomain.
For problems with spatially localized features, this construction has the advantage that local phenomena can often be represented in each subdomain using a relatively small number of local basis vectors.
Owing to this property, DD-G-ROMs have been increasingly applied to problems involving the evolution of localized physical phenomena, such as shock waves \citep{LUCIA2003917_snapshotPOD} and crack propagation \citep{KERFRIDEN2012154,KERFRIDEN2013169}, systems with repetitive geometrical components, as in component-based ROMs \citep{HUYNH2013213,SMETANA2016A3318,IOLLO2023115786,TADDEI2024113038,HUANG2022900064,IAPICHINO2016408,KAPTEYN20222986}, and multiphysics problems \citep{CORIGLIANO2013113,DISCACCIATI2023}.
However, even with DD-G-ROMs, the online phase can remain computationally demanding for complex phenomena.
Many existing studies do not parallelize the entire online phase and instead limit parallelization to the offline phase.
Although combinations of DD-G-ROMs with techniques such as the overlapping Schwarz method \citep{IOLLO2023115786}, the Schur complement method \citep{DECASTRO2024113282}, and the Lagrange multiplier method \citep{HOANG2021113997,DIAZ2024116943} suggest the potential for online-phase parallelization, the parallel efficiency achieved in practical computations has not yet been reported.
Moreover, these approaches are designed under the assumption that the domain decomposition used to construct local low-dimensional approximation spaces coincides with that used for parallel computing.
Under this assumption, increasing the number of local approximation-space subdomains requires access to high-performance computing resources capable of handling a correspondingly large number of parallel processes.
In addition, when the required reduced dimension varies among subdomains, the computational load can become imbalanced across parallel processes, potentially degrading the parallel efficiency.

\begin{figure}[t]
    \begin{center}
    \includegraphics[scale=0.5]{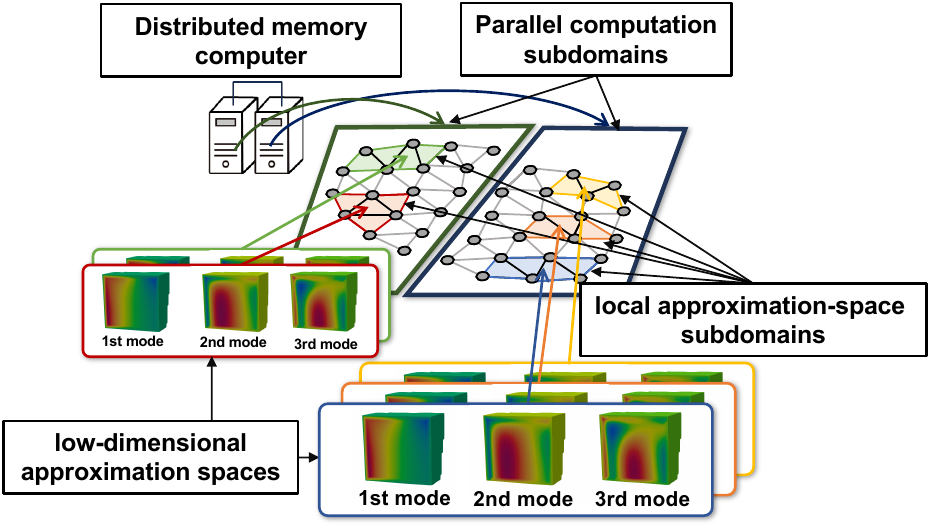}
    \end{center}
    \vspace{-10pt}
    \caption{\baselineskip=13pt Schematic overview of the proposed method.}
    \label{fig:intro}
    \vspace{-5pt}
\end{figure}

Against this background, this study aims to provide a general framework for alleviating memory constraints in large-scale DD-G-ROM analyses through distributed-memory parallelization of both the offline and online phases, while maintaining high parallel efficiency.
To this end, we propose a metagraph-based domain-decomposed Galerkin reduced-order model (MBDD-G-ROM).
The proposed framework introduces a novel graph-based representation of DD-G-ROMs and enables the subdomains used to construct local low-dimensional approximation spaces, referred to here as local approximation-space subdomains, to be specified flexibly and independently of the subdomains used for parallel computation.
We first propose a method for constructing graph representations corresponding to the discretization structures of the finite element node set and the local approximation-space subdomains.
This graph-based formulation provides a general description of DD-G-ROMs for arbitrary domain decompositions.
On this basis, we introduce a two-level graph-based domain-partitioning strategy that enables the local approximation-space subdomains to be specified independently of the subdomains used for distributed-memory parallelization.
Specifically, as illustrated in Fig.~\ref{fig:intro}, the finite element graph is first partitioned to define the local approximation-space subdomains.
A graph representing the connectivity among these subdomains, referred to as the metagraph, is then constructed and partitioned again to define the subdomains for distributed-memory parallelization.
In addition, by extending an overlapping domain-decomposition-based parallel computing method for finite element analyses \citep{YAGAWA1993495}, we propose a parallel computing method for dense block matrices that represent the interactions among the local approximation-space subdomains.
In this study, we adopt Local POD \citep{BAIGES201323}, in which low-dimensional approximation spaces are constructed by applying POD separately in each subdomain. Nevertheless, the proposed framework can be extended to other methods for constructing low-dimensional approximation spaces.

The remainder of this paper is organized as follows.
Section~2 describes the construction of ROMs based on snapshot POD.
Section~3 introduces the graph structures associated with DD-G-ROMs, which form an essential foundation of the proposed method, and presents a formulation applicable to arbitrary domain decompositions.
Section~4 introduces the concept of a metagraph based on these graph structures and proposes a method for model reduction and parallel computing based on a hierarchical graph structure. An example of load balancing based on weighted graphs is then presented.
Section~5 presents applications of the proposed method to an unsteady diffusion equation and the incompressible Navier--Stokes equations, and evaluates both ROM accuracy, including a quantitative assessment of the accuracy degradation caused by reduced-order modeling, and parallel computing performance.

%% file: include/equations.tex
\section{Snapshot POD-Based Galerkin Reduced-Order Model (ROM)}
This section describes Galerkin reduced-order models (G-ROMs), which compute solutions by projecting the governing equations onto a low-dimensional approximation space (a linear subspace or a nonlinear manifold), and the POD-Galerkin method, which is a G-ROM in which such a low-dimensional approximation space is constructed using proper orthogonal decomposition (POD).

\subsection{Full-Order Model}
This study considers partial differential equations that are first order in time. The spatial discretization is performed using the finite element method, and the temporal discretization is performed using an implicit time-integration scheme. 
The unknown vector $\bm u_{t+\Delta t}$ at time $t+\Delta t$ is determined by requiring the following weak residual equation to hold for any admissible test vector $\bm w$:
\begin{equation}\label{eq:residual}
  \bm w^{\mathrm{T}}\left(
     \bm M\,\dot{\bm u}_{t+\Delta t}
     + \bm f_{\mathrm{int}}\!\bigl(\bm u_{t+\Delta t}\bigr)-\bm f_{\mathrm{ext}}(t+\Delta t)\right)
     = 0.
\end{equation}
Here, $\Delta t$ denotes the time-step size, $\bm M\in\mathbb{R}^{d\,n_{\mathrm{FOM}}\times d\,n_{\mathrm{FOM}}}$ is the mass matrix, $\bm u_{t}\in\mathbb{R}^{d\,n_{\mathrm{FOM}}}$ is the unknown vector at time $t$, $\dot{\bm u}_{t}$ is the time derivative of $\bm u_{t}$, $n_{\mathrm{FOM}}$ is the number of nodes in the finite element discretization, and $d$ is the number of degrees of freedom per node. When an implicit time-integration scheme is adopted for the temporal discretization, the relation between $\bm u_{t}$ and $\dot{\bm u}_{t}$ is determined by the selected implicit time-integration scheme, such as the implicit Euler method or the generalized-$\alpha$ method \cite{JANSEN2000305}. In any of these schemes, however, the problem is reduced to solving Eq.~\eqref{eq:residual} at each time step. Furthermore, $\bm f_{\mathrm{int}}(\bm u_{t})\in\mathbb{R}^{d\,n_{\mathrm{FOM}}}$ is the internal force vector, $\bm f_{\mathrm{ext}}(t)\in\mathbb{R}^{d\,n_{\mathrm{FOM}}}$ is the external force vector acting at time $t$, and $\bm w\in\mathbb{R}^{d\,n_{\mathrm{FOM}}}$ is a vector derived from the test function. Although this study focuses on systems with a constant number of degrees of freedom per node, the formulation can be extended to systems in which the number of nodal degrees of freedom varies spatially, such as multiphysics problems represented by fluid--structure interaction, where the structural nodes have three degrees of freedom and the fluid nodes have four.

The nonlinear equation \eqref{eq:residual} is solved by Newton--Raphson iteration. For the unknown vector $\bm u_{t+\Delta t}$, the residual vector $\bm r \in \mathbb{R}^{d\,n_{\mathrm{FOM}}}$ and the Jacobian matrix $\bm J \in \mathbb{R}^{d\,n_{\mathrm{FOM}}\times d\,n_{\mathrm{FOM}}}$ are defined as
\begin{equation}\label{eq:Newton-raphson}
  \bm r(\bm u^{i}_{t+\Delta t})
  := 
       \bm M\,\dot{\bm u}^{i}_{t+\Delta t}
     + \bm f_{\mathrm{int}}\!\bigl(\bm u^{i}_{t+\Delta t}\bigr)
     - \bm f_{\mathrm{ext}}(t+\Delta t),
     \qquad \bm J(\bm u^{i}_{t+\Delta t})
:= \left.\frac{\partial \bm r(\bm u)}{\partial \bm u}\right|_{\bm u=\bm u^{i}_{t+\Delta t}},
\end{equation}
and the Newton--Raphson iteration is given by
\begin{align}\label{eq:NR_corrected}
  \bm w^{\mathrm{T}}\bm J (\bm u^{i}_{t+\Delta t})\,
  \Delta\bm u^{i}_{t+\Delta t}
  &= -\bm w^{\mathrm{T}}\bm r (\bm u^{i}_{t+\Delta t}),\\
  \bm u_{t+\Delta t}^{i+1}
  &= \bm u_{t+\Delta t}^{i}+\Delta\bm u^{i}_{t+\Delta t}. \label{eq:NR_ansvec}
\end{align}
Here, $i$ denotes the iteration index in the Newton--Raphson method. For a linear problem in which $\bm f_{\mathrm{int}}\!\bigl(\bm u_{t}\bigr) = \bm K \bm u_{t}$, it is sufficient to solve the system of linear equations obtained at each time step only once, without applying the Newton--Raphson method. Here, $\bm K\in \mathbb{R}^{d\,n_{\mathrm{FOM}}\times d\,n_{\mathrm{FOM}}}$ is the coefficient matrix in this linear representation.

\subsection{Galerkin ROM}
To express the high-dimensional unknown vector $\bm u\in\mathbb{R}^{d\,n_{\mathrm{FOM}}}$ in terms of the $k$-dimensional reduced coordinates $\bm q\in\mathbb{R}^{k}$, with $k\ll d\,n_{\mathrm{FOM}}$, we introduce a linear mapping $  \bm\Phi \in \mathbb{R}^{d\,n_{\mathrm{FOM}}\times k}$ and approximate the unknown vector as follows:
\begin{equation}\label{eq:mapping_linear}
  \bm u_t^i
  \approx
  \bm \Phi \bm q_t^i.
\end{equation}
The approximation in Eq.~\eqref{eq:mapping_linear} is used in combination with the Galerkin method, and the test vector is approximated as
\begin{equation}\label{eq:pod_approxi_weight}
\bm{w}\approx \bm\Phi \bm{\tilde{w}},
\end{equation}
where $\bm{\Tilde{w}}\in\mathbb{R}^{k}$ is an arbitrary vector.
Using the approximations in Eqs.~\eqref{eq:mapping_linear} and \eqref{eq:pod_approxi_weight}, Eqs.~\eqref{eq:NR_corrected} and \eqref{eq:NR_ansvec} are reduced to
\begin{align}\label{eq:NR_reduced}
   \tilde{\bm w}^{\mathrm{T}} \tilde{\bm J}(\bm q_{t+\Delta t}^{i})\,
  \Delta\bm q^{i}_{t+\Delta t}
  &= -\tilde{\bm w}^{\mathrm{T}}\tilde{\bm r}(\bm q_{t+\Delta t}^{i}),\\ \label{eq:NR_reduced_ansvec}
  \bm q_{t+\Delta t}^{i+1} &= \bm q_{t+\Delta t}^{i}+\Delta\bm q^{i}_{t+\Delta t},
\end{align}
where $\tilde{\bm J}(\bm q):= \bm \Phi^{\mathrm{T}} \bm J(\bm \Phi \bm q) \bm \Phi$ and $\tilde{\bm r}(\bm q):=\bm \Phi^{\mathrm{T}} \bm r (\bm \Phi \bm q)$. In Eqs.~\eqref{eq:NR_reduced} and \eqref{eq:NR_reduced_ansvec}, the unknown vector is $\Delta \bm q\in\mathbb{R}^{k}$. Thus, the dimension of the system is reduced from $d\,n_{\mathrm{FOM}}$ to $k$, enabling fast solution.

Alternatively, nonlinear-manifold ROMs introduce a nonlinear mapping $\bm\Gamma : \mathbb{R}^{k}\to\mathbb{R}^{d\,n_{\mathrm{FOM}}}$ and approximate the unknown vector $\bm u(t)$ as $\bm u(t) \approx \bm\Gamma\bigl(\bm q(t)\bigr)$ to reduce the governing equations. Recent studies have actively investigated approaches based on quadratic manifolds \citep{RUTZMOSER2017196,BARNETT2022111348,GEELEN2023115717} and autoencoders \cite{KIM2024116978,LEE2020108973}. However, because the primary objective of this study is to propose a new parallel analysis framework based on domain decomposition and metagraphs, we adopt the linear-subspace construction method based on POD, introduced in Section~\ref{syou:POD-G}, for simplicity of discussion and ease of verification.

\subsection{POD-Galerkin Method}\label{syou:POD-G}
In this study, we construct and solve ROMs using the POD-Galerkin method based on Eqs.~\eqref{eq:mapping_linear} and \eqref{eq:pod_approxi_weight}, where the basis is extracted in a data-driven manner using snapshot POD.

First, the snapshot matrix $\ \bm{S}\in \mathbb{R}^{d\,n_{\mathrm{FOM}} \times n_{\mathrm{snap}}}$ is formed from $n_{\mathrm{snap}}$ snapshot data vectors $\bm{d}_{(i)}\in \mathbb{R}^{d\,n_{\mathrm{FOM}}}\ (i = 1,2,\ldots,n_{\mathrm{snap}})$ as follows:
\begin{equation}
    \bm{S}=\left[\bm{d}_{(1)},\bm{d}_{(2)},\ldots,\bm{d}_{(n_{\mathrm{snap}})}\right].
\end{equation}
We then consider obtaining ${n_{\mathrm{POD}}}$ basis vectors $\bm{\Phi}_i\ \in \mathbb{R}^{d\,n_{\mathrm{FOM}}}\ (i=1,2,\ldots,n_{\mathrm{POD}})$ that represent the features of the snapshot matrix $\ \bm{S}$. POD determines the basis of the subspace that minimizes the error between each snapshot and its orthogonal projection onto the subspace, as follows:
\begin{equation}\label{eq:pod_equ}
    \{\bm{\Phi}_i\}_{i=1}^{n_{\mathrm{POD}}}=\underset{\bm{g}_{i}\ (i=1,2,\ldots,n_{\mathrm{POD}})}{\mathrm{arg\ min}}{\sum_{j=1}^{n_{\mathrm{snap}}}\left|\left|\bm{d}_{(j)}
    -\sum_{l=1}^{{n_{\mathrm{POD}}}}{\bm{g}}_l\bm{g}_l^{\mathrm{T}}\bm{d}_{(j)}\right|\right|^2},
\end{equation}
\begin{equation}
    \mathrm{with} \ \  \bm{G}^{\mathrm{T}}\bm{G}=\bm{I},\ \ \bm{G}=[\bm{g}_1,\bm{g}_2,\ldots,\bm{g}_{n_{\mathrm{POD}}}]\in \mathbb{R}^{d\,n_{\mathrm{FOM}} \times {n_{\mathrm{POD}}}}.\notag
\end{equation}
The singular value decomposition of the snapshot matrix $\bm{S}$ is written as
\begin{equation}\label{eq:SVD}
    \bm{S}=\Hat{\bm{\Phi}}\bm{\varSigma} \bm{W}^{\mathrm{T}}.
\end{equation}
In Eq.~\eqref{eq:SVD}, each column of $\Hat{\bm{\Phi}}=[\bm{\Phi}_1,\bm{\Phi}_2,\ldots,\bm{\Phi}_r]\in \mathbb{R}^{d\,n_{\mathrm{FOM}} \times r}\ (r = \mathrm{\mathrm{rank}}\ \bm{S})$ is a left singular vector, and $\bm{\varSigma}=\mathrm{diag}{(\sigma_{i})}\ \in \mathbb{R}^{r\times r}\ (i=1,2,\ldots,r)$ is a matrix whose diagonal entries are the singular values $\sigma_i$. The quantity $\sigma_i$ denotes the $i$th singular value, and the singular values are nonnegative and ordered in descending magnitude with respect to the index $i$. Each column of $\bm{W}=[\bm{w}_1,\bm{w}_2,\ldots,\bm{w}_r]\in \mathbb{R}^{n_{\mathrm{snap}} \times r}$ is a right singular vector. The first ${n_{\mathrm{POD}}}$ columns of $\Hat{\bm{\Phi}}$ corresponding to the largest singular values are selected to obtain the POD basis $  \bm\Phi
  = \bigl[\bm{\Phi}_{1},\bm{\Phi}_2,\dots, \bm{\Phi}_{n_{\mathrm{POD}}}\bigr]\in\mathbb{R}^{d\,n_{\mathrm{FOM}}\times {n_{\mathrm{POD}}}}$. Because the singular value $\sigma_i$ indicates the contribution of the $i$th basis vector, the number of basis vectors ${n_{\mathrm{POD}}}$ is determined using an appropriate threshold $\varepsilon_{\mathrm{POD}}$ so as to satisfy
\begin{equation}\label{eq:base_sel}
    \frac{\sum_{i=1}^{n_{\mathrm{POD}}}{\sigma_i}}{\sum_{i=1}^{r}{\sigma_i}} >
    1-\varepsilon_{\mathrm{POD}}.
\end{equation}
In the POD-Galerkin method, model reduction is performed using the POD basis $\bm \Phi$ obtained in this manner, according to Eqs.~\eqref{eq:mapping_linear} and \eqref{eq:pod_approxi_weight}.

\section{Graph-Based Generalization of Domain-Decomposed ROM}\label{syou:Graph-BasedROM}
Existing studies on DD-G-ROMs (domain-decomposed G-ROMs) have proposed frameworks in which a physical domain or algebraic degrees of freedom are partitioned into multiple subdomains, and a low-dimensional approximation space is constructed for each subdomain. Representative approaches include methods that construct a low-dimensional approximation space for each subdomain and couple the subdomains continuously based on the finite element method \citep{BAIGES201323}, methods that construct a low-dimensional approximation space for each component \citep{HUYNH2013213}, algebraic non-overlapping domain decomposition \citep{HOANG2021113997}, and optimization-based approaches \citep{Prusak2023,TADDEI2024113038}. These studies have demonstrated that DD-G-ROMs are effective in terms of representing local solution structures, localizing training data in the offline phase, and extending ROMs to large-scale problems. 

However, many of these approaches focus on predefined components, algebraic domain decomposition, or optimization-based coupling of subdomain interfaces, and relatively few studies have directly connected arbitrary mesh partitioning to the offline- and online-phase construction procedures of DD-G-ROMs. A closely related method is the domain-decomposition non-intrusive reduced-order model (DDNIROM) \citep{Xiao2019CF}. In DDNIROM, the computational domain is divided into multiple subdomains, low-dimensional approximation spaces based on POD are constructed on individual subdomains, and local approximation functions in the reduced spaces are constructed using Gaussian process regression, thereby representing localized nonlinear flow features. In particular, DD-ROMs that are not limited to simple geometric decompositions have been realized by domain decompositions that account for flow-field features through nodal weights \citep{Xiao2019CF}. However, DDNIROM belongs to the class of non-intrusive ROMs in which the online solution is predicted using machine-learning models.
Its formulation therefore differs from the DD-G-ROM considered in this study, in which the discrete FOM equations are projected onto low-dimensional approximation spaces and the reduced solution over the entire computational domain is reconstructed based on the governing equations. In this section, we introduce graph concepts corresponding to DD-G-ROMs and to the sparsity structures of their matrices, and present a formulation of DD-G-ROMs generalized to arbitrary domain decompositions. We then detail the formulation for the case in which Local POD is used.

We begin by describing the graphs considered in this study. A graph $G=(V,E)$ consists of a node set $V=\{v_1,v_2,\ldots,v_n\}$ and an edge set $E$, where $n$ denotes the number of nodes $|V|$. When nodes $v_i$ and $v_j$ are connected, the connection is denoted by $(v_i,v_j)$ and is referred to as an edge. The graphs defined in this paper are assumed to include self-edges $(v_i,v_i)\ (i=1,2,\ldots,|V|)$ and are treated as undirected graphs. In addition, the set $\mathcal{V}$ consisting of the indices of the nodes in $V$ is defined as
\begin{equation}\label{eq:index_set}
\mathcal{V} := \bigl\{\, i\mid v_i \in V \,\bigr\}.
\end{equation}

\subsection{Graph-Based Domain Decomposition Method}\label{syou:DD}

\begin{figure*}[t]
    \begin{center}
    \includegraphics[scale=0.6]{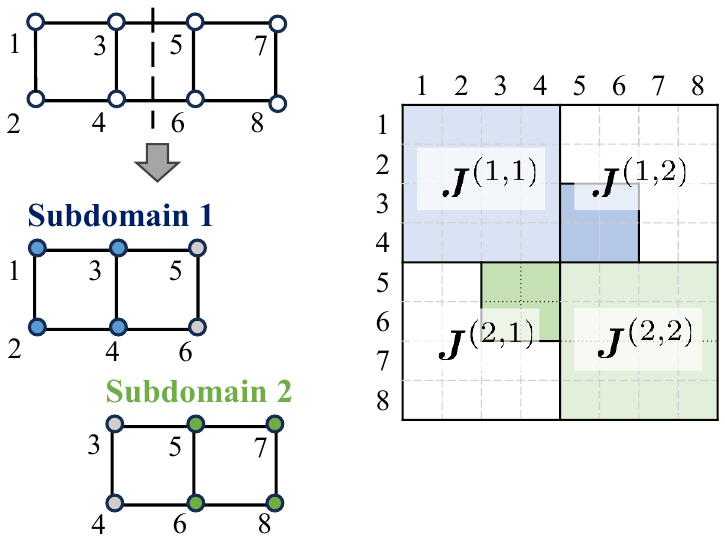}
    \end{center}
    \vspace{-10pt}
    \caption{\baselineskip=13pt Relationship between subdomains and the block-matrix structure.}
    \label{fig:DDM_matrix}
\end{figure*}

A calculation-point graph is a graph structure that represents interactions among computational points and corresponds to the nonzero structure of the coefficient or Jacobian matrix appearing in Eq.~\eqref{eq:NR_corrected} under consideration. For example, a calculation-point graph whose nodes are the nodes of a finite element mesh corresponds to the nonzero structure of the system of linear equations obtained by finite element discretization. Using such a graph structure, domain decomposition can be defined for an arbitrary node set with a given adjacency relation.

First, under the following definitions, the node set $V$ of the graph is partitioned without overlap into $N$ node sets $V^{(i)}\ (i = 1,2,\ldots,N)$:
\begin{align}
    V^{(i)}\cap V^{(j)}&=\emptyset\quad (i\neq j),\label{metis2}\\
    \underset{i = 1}{\overset{N}{\bigcup}} V^{(i)} &= V. \label{metis3}
\end{align}
The set of graph nodes belonging to each subdomain is referred to as the internal node set, and the number of internal nodes belonging to subdomain $i$ is denoted by $|V^{(i)}|$. The node set $V^{(i)}$ obtained by domain decomposition is referred to as subdomain $i$. The set $\mathcal{V}^{(i)}$ consisting of the indices of the nodes in $V^{(i)}$ is defined as
\begin{equation}
\mathcal V^{(i)} := \{\, j \in \mathcal V \mid v_j \in V^{(i)} \,\}.
\end{equation}
The node set $\bar{V}^{(i)}$ that accounts for the overlap region is given by
\begin{equation}\label{metis_overlap1}
    \bar{V}^{(i)}=V^{(i)}\cup V^{(i)}_{\mathrm{ovl}},
\end{equation}
where $V^{(i)}_{\mathrm{ovl}}$ is defined as
\begin{equation}
V_{\mathrm{ovl}}^{(i)} := \left\{\, v_j \mid v_j \notin V^{(i)},\ \exists v_k \in V^{(i)}\ \text{s.t.}\ (v_k, v_j) \in E \,\right\}.
\end{equation}
Subdomain $i$ is said to be adjacent to subdomain $j$ when the following condition holds:
\begin{equation}\label{eq:grapheq1}
V^{(i)}\cap V^{(j)}_{\mathrm{ovl}} \neq \emptyset.
\end{equation}
The left side of Fig.~\ref{fig:DDM_matrix} shows a schematic of a finite element mesh corresponding to the node set $V$ when it is partitioned into two subdomains.

\begin{figure*}[t]
    \begin{center}
    \includegraphics[scale=0.6]{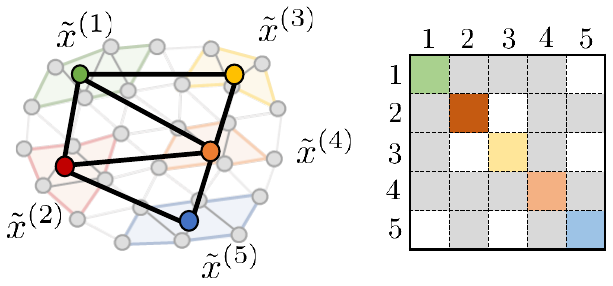}
    \end{center}
    \vspace{-10pt}
    \caption{\baselineskip=13pt Relationship between the metagraph and the block-matrix structure (nonzero blocks are shown in color or gray).}
    \label{fig:metagraph}
\end{figure*}

Next, we define a metagraph that corresponds to the adjacency relation among subdomains obtained by graph partitioning of a calculation-point graph. As shown on the left side of Fig.~\ref{fig:metagraph}, a new metanode $\Tilde{v}_i$ is defined for each internal node set $V^{(i)}$ obtained by partitioning the calculation-point graph. For these metanodes, a metagraph $\Tilde{G}=(\Tilde{V},\Tilde{E})$ representing the adjacency relation among the subdomains can be defined. The node set $\Tilde{V}=\{\Tilde{v}_1,\Tilde{v}_2,\ldots,\Tilde{v}_N\}$ represents the partitioned subdomains, and the edge set $\Tilde{E} = \bigl\{(\Tilde{v}_{i}, \Tilde{v}_{j})\ |\ \Tilde{v}_i,\Tilde{v}_j\in \Tilde{V},\, i=j
\ \lor\ V^{(i)}\cap V^{(j)}_{\mathrm{ovl}} \neq \emptyset\bigr\}$ represents the self-edges corresponding to the subdomains themselves and the adjacency relations among the subdomains.

We next clarify the correspondence between a matrix $\bm{A}$ and a vector $\bm{x}$ before and after graph partitioning based on the calculation-point graph. First, consider a vector $\bm{x}$ of size $\sum_{i\in\mathcal{V}}d_{i}$ for the calculation-point graph $G=(V,E)$:
\begin{equation}
\bm{x}
=
\begin{bmatrix}
\bm{x}_1 \\
\bm{x}_2 \\
\vdots \\
\bm{x}_{|V|}
\end{bmatrix}
\in\mathbb{R}^{\sum_{i\in\mathcal{V}}d_{i}}.
\end{equation}
Here, $d_i$ denotes the number of degrees of freedom at computational point $v_i$, and $\bm x_i$ is the vector at computational point $v_i$. A square matrix $\bm{A}$ of size $\sum_{i\in\mathcal{V}}d_{i} \times \sum_{i\in\mathcal{V}}d_{i}$ associated with the calculation-point graph $G=(V,E)$ is expressed as the following block matrix:
\begin{equation}
\bm{A}
=
\begin{bmatrix}
\bm{A}_{1,1} & \bm{A}_{1,2} & \cdots & \bm{A}_{1,|V|} \\
\bm{A}_{2,1} & \bm{A}_{2,2} & \cdots & \bm{A}_{2,|V|} \\
\vdots & \vdots & \ddots & \vdots \\
\bm{A}_{|V|,1} & \bm{A}_{|V|,2} & \cdots & \bm{A}_{|V|,|V|}
\end{bmatrix}
\in \mathbb{R}^{\sum_{i\in\mathcal{V}}d_{i} \times \sum_{j\in\mathcal{V}}d_{j}}.
\end{equation}
The block $\bm{A}_{i,j}\in\mathbb{R}^{d_i\times d_j}$ represents the action from computational point $v_j$ to computational point $v_i$, and $\bm{A}_{i,j}$ is a zero matrix when $(v_i,v_j)\notin E$.

Partitioning the calculation-point graph $G=(V,E)$ yields a partitioned graph $G^{(i)}$ corresponding to subdomain $i$. The vector $\bm{x}^{(i)}$ corresponding to the internal node set $V^{(i)}$ of subdomain $i$ is expressed as
\begin{equation}
\bm{x}^{(i)} =
\begin{bmatrix}
\bm{x}^{(i)}_{1} \\
\bm{x}^{(i)}_{2} \\
\vdots \\
\bm{x}^{(i)}_{|V^{(i)}|}
\end{bmatrix}
\in \mathbb{R}^{\sum_{j\in\mathcal{V}^{(i)}}d_{j}},
\end{equation}
where $\bm x^{(i)}_p$ is the vector corresponding to computational point $v_{(\mathrm{idx}(i,p))}$ in subdomain $i$. The function $\mathrm{idx}(i,a_\mathrm{L})$ maps the local number $a_\mathrm{L}$ in subdomain $i$ to the corresponding global number $a_\mathrm{G}$, that is, $a_\mathrm{G}=\mathrm{idx}(i,a_\mathrm{L})$. In addition, the matrix $\bm{A}^{(i,j)}$, which represents the action of the internal node set $V^{(j)}$ of subdomain $j$ on the internal node set $V^{(i)}$ of subdomain $i$, is expressed as the following block matrix:
\begin{equation}
\bm{A}^{(i,j)}
=
\begin{bmatrix}
\bm{A}^{(i,j)}_{1,1} & \bm{A}^{(i,j)}_{1,2} & \cdots & \bm{A}^{(i,j)}_{1,|V^{(j)}|} \\
\bm{A}^{(i,j)}_{2,1} & \bm{A}^{(i,j)}_{2,2} & \cdots & \bm{A}^{(i,j)}_{2,|V^{(j)}|} \\
\vdots & \vdots & \ddots & \vdots \\
\bm{A}^{(i,j)}_{|V^{(i)}|,1} & \bm{A}^{(i,j)}_{|V^{(i)}|,2} & \cdots & \bm{A}^{(i,j)}_{|V^{(i)}|,|V^{(j)}|}
\end{bmatrix}\in \mathbb{R}^{\sum_{k\in\mathcal{V}^{(i)}}d_{k} \times \sum_{l\in\mathcal{V}^{(j)}}d_{l}}.
\end{equation}
Here, $\bm{A}^{(i,j)}_{p,q}$ represents the action from the $q$th computational point $v_{(\mathrm{idx}(j,q))}\in V^{(j)}$ in $V^{(j)}$ to the $p$th computational point $v_{(\mathrm{idx}(i,p))}\in V^{(i)}$ in $V^{(i)}$, and $\bm{A}^{(i,j)}_{p,q}$ is a zero matrix when $(v_{(\mathrm{idx}(i,p))},v_{(\mathrm{idx}(j,q))})\notin E$.

The node number before domain decomposition is defined as the global number
$a_\mathrm{G}$, and the node number in each subdomain after domain decomposition
is defined as the local number $a_\mathrm{L}$. The relationships between the
matrices $\bm{A}$ and $\bm{A}^{(i,j)}$ and between the vectors $\bm{x}$ and
$\bm{x}^{(i)}$ before and after domain decomposition are expressed using the
matrix
$\bm{P}^{(i)}\in
\mathbb{R}^{\sum_{j\in\mathcal{V}}d_{j}
\times
\sum_{k\in\mathcal{V}^{(i)}}d_{k}}$.
The matrix $\bm{P}^{(i)}$ is a local-to-global Boolean assembly matrix that
embeds a local vector defined on subdomain $i$ into the corresponding entries
of the global vector. Its transpose ${\bm{P}^{(i)}}^{\mathrm{T}}$ restricts a
global vector to subdomain $i$. These relationships are written as follows:
\begin{equation}\label{eq:para_eq_mat}
\bm{A}
=
\sum_{i=1}^{N}\sum_{j=1}^{N}
\bm{P}^{(i)}
\bm{A}^{(i,j)}
{\bm{P}^{(j)}}^{\mathrm{T}},
\end{equation}
\begin{equation}\label{eq:para_eq_vec}
    \bm{x}=\sum_{i=1}^{N}\bm{P}^{(i)}\bm{x}^{(i)}.
\end{equation}
When subdomains $i$ and $j$ are adjacent, the block matrix $\bm{A}^{(i,j)}$ representing the interaction between them is nonzero. Because the metagraph represents the adjacency relation among subdomains, the nonzero structure of the block matrices $\bm{A}^{(i,j)}$ of the matrix $\bm{A}$ is uniquely determined by referring to the metagraph, as illustrated in Fig.~\ref{fig:metagraph}. The matrix entries of $\bm{P}^{(i)}$ are given by
\begin{equation}\label{eq:permutation_mat}
\bm{P}^{(i)}_{m,\,n}=
\begin{cases}
\bm I \in \mathbb{R}^{d_{m}\times d_{m}} 
& m=\mathrm{idx}\!\left(i, n\right)\\[2pt]
\bm 0 \in \mathbb{R}^{d_{m}\times d_{\mathrm{idx}(i, n)}}
& \text{otherwise}
\end{cases},
\end{equation}
where $d_m$ is the number of degrees of freedom at computational point $v_m$.

As a concrete example of a calculation-point graph, consider a graph whose nodes are the nodes in the finite element method and whose edges represent their adjacency relations. Let $V\ (|V| = n_{\mathrm{FOM}})$ be the finite element node set, and let $N_a$ and $N_b$ be the basis functions corresponding to nodes $a$ and $b$, respectively. Then, a pair of finite element nodes $( v_a,v_b )$ for which the product of these basis functions is nonzero on the same element is regarded as an edge. The edge set $E$ derived from interactions among finite elements is then defined as
\begin{equation}
  E=\bigl\{( v_a,v_b )
         \mid \,
         \exists\bm x\in \mathbb{R}^3:\,
         N_a(\bm x)\,N_b(\bm x)\neq0
   \bigr\}.
\end{equation}
For the finite element graph $G = (V,E)$, subgraphs are obtained by graph partitioning based on Eqs.~\eqref{metis2} and \eqref{metis3}. A metagraph $\Tilde{G}=(\Tilde{V},\Tilde{E})$ can also be defined.

The Jacobian $\bm{J}\in \mathbb{R}^{d\,n_{\mathrm{FOM}}\times d\,n_{\mathrm{FOM}}}$ corresponding to the node set $V$ before domain decomposition, defined in Eq.~\eqref{eq:Newton-raphson}, is expressed using the Jacobian $\bm{J}^{(i,j)}\in \mathbb{R}^{d\,n_{\mathrm{FOM}}^{(i)} \times d\,n_{\mathrm{FOM}}^{(j)}}$ corresponding to the internal node sets $V^{(i)}$ and $V^{(j)}$ of subdomains $i$ and $j$ $(|V^{(i)}|=n_{\mathrm{FOM}}^{(i)},|V^{(j)}|=n_{\mathrm{FOM}}^{(j)})$, as follows for a decomposition into $N$ subdomains:
\begin{equation}\label{eq:L-POD_block}
    \bm{J}=\left[
    \begin{array}{cccc}
      \bm{J}^{(1,1)}&\bm{J}^{(1,2)}&\cdots&\bm{J}^{(1,N)}\\
        \bm{J}^{(2,1)}&\bm{J}^{(2,2)}&\cdots&\bm{J}^{(2,N)}\\
                \vdots&\vdots&\ddots&\vdots\\
        \bm{J}^{(N,1)}&\bm{J}^{(N,2)}&\cdots&\bm{J}^{(N,N)}\\
    \end{array}
  \right].
\end{equation}
Here, $\bm{J}^{(i,j)}$ is the matrix representing the action received by metanode $i$ from metanode $j$, and $\bm{J}^{(i,j)}$ is a zero matrix when $(\tilde{v}_i,\tilde{v}_j)\notin \tilde{E}$.

The residual vector $\bm{r}\in \mathbb{R}^{d\,n_{\mathrm{FOM}}}$ corresponding to the node set $V$ before domain decomposition is expressed using the residual vector $\bm{r}^{(i)}\in \mathbb{R}^{d\,n_{\mathrm{FOM}}^{(i)}}$ corresponding to the internal node set $V^{(i)}$ of subdomain $i$ as
\begin{equation}\label{eq:DD_vec}
  \bm{r}=
\left[
    \begin{array}{cccc}
      \bm{r}^{(1)}\\
      \bm{r}^{(2)}\\
      \vdots\\
      \bm{r}^{(N)}\\
    \end{array}
  \right].
\end{equation}

\subsection{Graph-Based Generalization of Local POD (L-POD)}\label{syou:Local POD}

L-POD \cite{BAIGES201323} performs singular value decomposition separately for the snapshot matrix associated with each subdomain and assembles the resulting local bases into a global, block-structured reduced basis. The governing equations of the FOM are then projected onto the approximation space spanned by this basis using Galerkin projection. Because the singular value decompositions are performed on smaller, subdomain-level snapshot matrices, L-POD can reduce the offline cost of basis construction and may permit a larger total number of basis vectors to be retained than global POD. An increase in the reduced dimension may, however, increase the computational cost of the online ROM, while potentially improving its approximation accuracy. Moreover, the local basis dimension can be selected independently for each subdomain, allowing the reduced approximation space to be tailored to spatial variations in solution complexity.

\begin{figure*}[t]
    \begin{center}
    \includegraphics[scale=0.6]{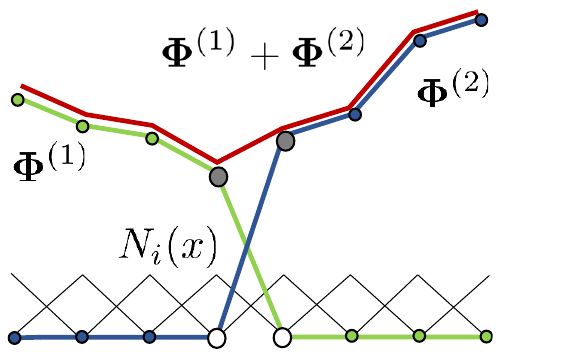}
    \end{center}
    \vspace{-10pt}
    \caption{\baselineskip=13pt Schematic illustration of local POD basis functions and their overlap across adjacent subdomains.}
    \label{fig:LPOD}
\end{figure*}

First, the singular value decomposition in Eq.~\eqref{eq:SVD} is applied to each local snapshot matrix $\bm{S}^{(i)}\in \mathbb{R}^{d\,n_{\mathrm{FOM}}^{(i)}\times n_{\mathrm{snap}}}\ (i = 1,2,\ldots,N_{\mathrm{POD}})$ corresponding to the internal node set $V^{(i)}$ of subdomain $i$:
\begin{equation}\label{eq:para_SVD}
    \bm{S}^{(i)}=\Hat{\bm{\Phi}}^{(i)}\bm{\Sigma}^{(i)}{\bm{W}^{(i)}}^{\mathrm{T}}.
\end{equation}
Here, $N_{\mathrm{POD}}$ denotes the number of subdomains used to obtain POD bases. Among the columns of the matrix $\Hat{\bm{\Phi}}^{(i)}\in \mathbb{R}^{d\,n_{\mathrm{FOM}}^{(i)}\times r^{(i)}}$ of left singular vectors, where $r^{{(i)}} = \mathrm{rank}\ \bm{S}^{(i)}$, the first $n_{\mathrm{POD}}^{(i)}$ columns are selected and used as the POD basis $\bm{\Phi}^{(i)}\in \mathbb{R}^{d\,n_{\mathrm{FOM}}^{(i)}\times n_{\mathrm{POD}}^{(i)}}$. Because the bases are obtained independently in the individual subdomains, the POD basis $\bm{\Phi}$ over the entire domain is given by
\begin{equation}\label{eq:localPODmodes}
    \bm{\Phi}= \left[
    \begin{array}{cccc}
    \bm{\Phi}^{(1)}  &0&\cdots&0\\
    0&\bm{\Phi}^{(2)}  &\cdots&0\\
     \vdots&\vdots& \ddots &   \vdots\\
    0&0&\cdots&\bm{\Phi}^{({N_{\mathrm{POD}}})} \\
    \end{array}
  \right].
\end{equation}
Figure~\ref{fig:LPOD} shows an example of L-POD for a one-dimensional problem with linear finite element basis functions. An L-POD basis is a linear combination of finite element basis functions associated with the nodes in subdomain $i$. Therefore, each row component of $\bm \Phi^{(i)}$ corresponds to a nodal value on the finite element mesh. In Fig.~\ref{fig:LPOD}, the green POD basis belongs to the left subdomain, and the blue POD basis belongs to the right subdomain. The red line represents the sum of the green and blue POD bases. Because the original two local POD bases are continuous, the red global POD basis is also continuous. In addition, there exists an overlap region in which both the left and right local POD bases are nonzero.

The coefficient vector $\bm{q}$ over the entire domain is expressed using the coefficient vectors $\bm{q}^{(i)}\in \mathbb{R}^{n_{\mathrm{POD}}^{(i)}}$ corresponding to the partitioned subdomains as
\begin{equation}\label{eq:reduced_ans}
\bm{q}=
\left[
    \begin{array}{cccc}
      \bm{q}^{(1)}\\
      \bm{q}^{(2)}\\
      \vdots\\
      \bm{q}^{({N_{\mathrm{POD}}})}\\
    \end{array}
  \right].
\end{equation}
In L-POD, as in conventional POD, the solution vector is approximated using this $\bm{q}$ and $\bm{\Phi}$. The test function is also approximated using $\bm{\Phi}$ by applying the Galerkin method, yielding the reduced systems of linear equations in Eq.~\eqref{eq:NR_reduced}.

The reduced Jacobian $ \tilde{\bm{J}} = \bm{{\Phi}}^{\mathrm{T}}\bm{J}\bm{{\Phi}}$ has the following block matrix structure:
\begin{equation}\label{eq:Local POD1}
     \tilde{\bm{J}}=\left[
    \begin{array}{cccc}
     \tilde{\bm{J}}^{(1,1)}& \tilde{\bm{J}}^{(1,2)}&\cdots& \tilde{\bm{J}}^{(1,{N_{\mathrm{POD}}})}\\
       \tilde{\bm{J}}^{(2,1)}& \tilde{\bm{J}}^{(2,2)}&\cdots& \tilde{\bm{J}}^{(2,{N_{\mathrm{POD}}})}\\
                \vdots&\vdots&\ddots&\vdots\\
       \tilde{\bm{J}}^{({N_{\mathrm{POD}}},1)}& \tilde{\bm{J}}^{({N_{\mathrm{POD}}},2)}&\cdots& \tilde{\bm{J}}^{({N_{\mathrm{POD}}},{N_{\mathrm{POD}}})}\\
    \end{array}
  \right].
\end{equation}
Each block matrix $\tilde{\bm{J}}^{(i,j)} \in \mathbb{R}^{n_{\mathrm{POD}}^{(i)}\times n_{\mathrm{POD}}^{(j)}}$ of the reduced Jacobian is expressed as a product with the POD bases corresponding to the associated nodes and is computed as follows:
\begin{equation}\label{eq:Local_POD_eq1}
   \tilde{\bm{J}}^{(i,j)}={\bm{\Phi}^{(i)}}^{\mathrm{T}}\bm{J}^{(i,j)}\bm{\Phi}^{(j)}.
\end{equation}
When subdomains $i$ and $j$ are adjacent, the matrix product ${\bm{\Phi}^{(i)}}^{\mathrm{T}}\bm{J}^{(i,j)}\bm{\Phi}^{(j)}$ becomes nonzero; therefore, Eq.~\eqref{eq:Local POD1} has a block nonzero structure similar to that of Eq.~\eqref{eq:L-POD_block}.

The reduced residual vector $\tilde{\bm{r}} = \bm{\Phi}^{\mathrm{T}}\bm{{r}}$ is computed as the product with the POD bases belonging to the corresponding subdomains:
\begin{equation}\label{eq:reduced_rhs}
\tilde{\bm{r}}=
\left[
    \begin{array}{ccc}
      \tilde{\bm{r}}^{(1)}\\
            \tilde{\bm{r}}^{(2)}\\
      \vdots\\
      \tilde{\bm{r}}^{({N_{\mathrm{POD}}})}\\
    \end{array}
  \right]\\
  =
\left[
    \begin{array}{ccc}
      {\bm{\Phi}^{(1)}}^{\mathrm{T}}\bm{r}^{(1)} \\
            {\bm{\Phi}^{(2)}}^{\mathrm{T}}\bm{r}^{(2)} \\
      \vdots\\
      {\bm{\Phi}^{({N_{\mathrm{POD}}})}}^{\mathrm{T}}\bm{r}^{({N_{\mathrm{POD}}})} \\
    \end{array}
  \right].
\end{equation}
Here, $\tilde{\bm{r}}^{(i)} \in \mathbb{R}^{n_{\mathrm{POD}}^{(i)}}$ is the reduced residual vector in subdomain $i$, and $\bm{r}^{(i)} \in \mathbb{R}^{d\,n_{\mathrm{FOM}}^{(i)}}$ is the residual vector in subdomain $i$.

\section{Metagraph-Based Domain-Decomposed Galerkin ROM}
As outlined in the Introduction, distributed-memory parallelization of L-POD over arbitrary domain decompositions has not yet been systematically formulated. Because L-POD has a matrix structure in which the block matrices corresponding to adjacent POD computation subdomains, namely the subdomains used to compute POD bases, are nonzero, information transfer among POD computation subdomains must be considered. In addition, in distributed-memory parallel computing, information transfer among the subdomains used for distributed-memory parallelization must also be considered for data communication among parallel processes. Therefore, distributed-memory parallel computation of L-POD requires handling information transfer over two types of hierarchical subdomains, and the key issue is how to realize this mathematically. Section \ref{syou:Graph-BasedROM} presented a unified graph-based formulation for DD-G-ROMs defined over arbitrary domain decompositions by introducing the concepts of calculation-point graphs and metagraphs. Building on this formulation, we propose the metagraph-based DD-G-ROM (MBDD-G-ROM) and formulate a hierarchical graph-partitioning method that decouples the POD computation subdomains from the parallel computation subdomains, thereby enabling distributed-memory parallelization of L-POD over arbitrary domain decompositions.

\subsection{Metagraph-Based Parallel Computation}\label{syou:parallel_comp}

Methods for solving systems of linear equations can be broadly classified into direct and iterative methods. Direct methods typically rely on factorization or elimination steps that are less straightforward to parallelize efficiently on distributed-memory systems, whereas iterative methods are readily expressed in terms of parallelizable algebraic kernels such as matrix--vector products and inner products. This study targets large-scale problems and therefore adopts iterative methods with high parallel efficiency for parallel computers. From the algorithmic structure of iterative methods, their parallel computation can be realized by implementing parallel vector sums, vector inner products, and matrix--vector products. Except for the data communication required for parallel computation, including global reductions for inner products, these parallel algebraic operations can be performed independently in each subdomain using the block matrices and vectors assigned to that subdomain, and thus they can achieve high parallel efficiency. This subsection describes parallel computation procedures for these linear algebra operations based on the calculation-point graph. In this study, parallel computation is performed using the Message Passing Interface (MPI) \cite{MPI}. MPI is a communication standard for parallel computing and enables highly scalable parallel implementations.

\begin{figure*}[t]
    \begin{center}
    \includegraphics[scale=0.42]{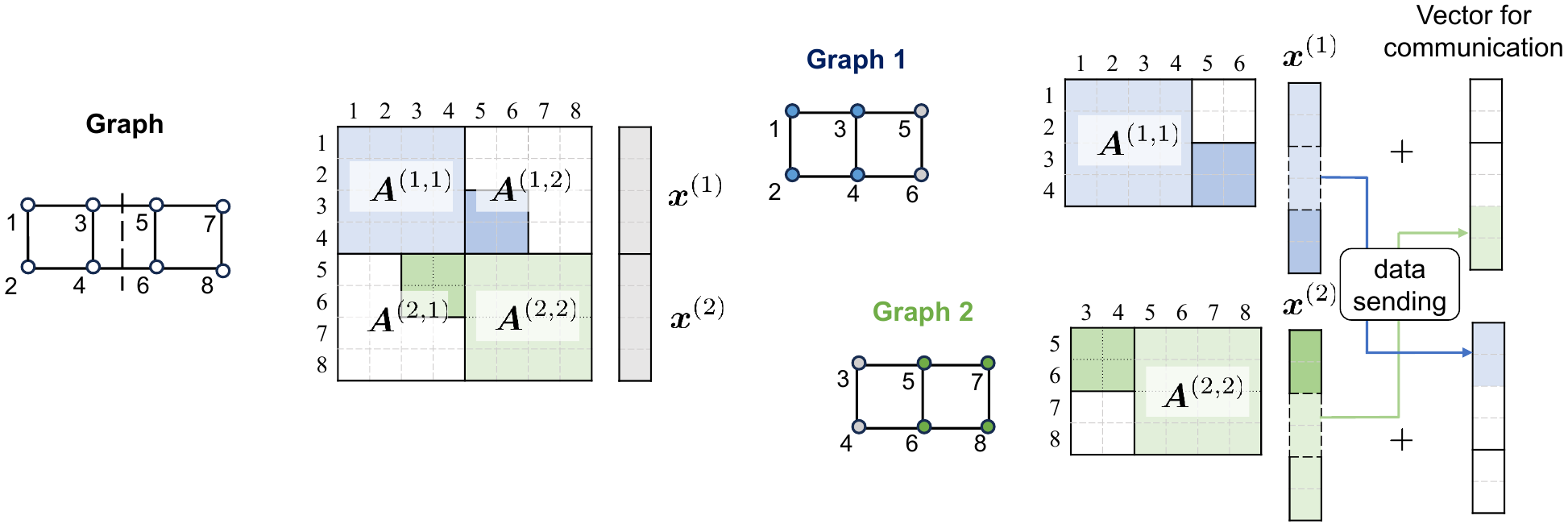}
    \end{center}
    \vspace{-10pt}
    \caption{\baselineskip=13pt  Schematic illustration of parallel sparse matrix--vector multiplication on graph $G$.}
    \label{fig:matvec}
\end{figure*}

When the node set $V$ of the graph is partitioned into $N_{\mathrm{par}}$ node sets $V^{(i)}$, a parallel process is assigned to the index $i$ representing the $i$th subdomain in Eqs.~(\ref{eq:para_eq_mat}) and (\ref{eq:para_eq_vec}). The parallel process associated with index $i$ is assigned the subgraph $G^{(i)}$ and stores the block matrices $\bm{A}^{(i,j)}$ and the vector $\bm{x}^{(i)}$. Here, $N_{\mathrm{par}}$ denotes the number of subdomains to which parallel processes are assigned.

The vector sum $\bm{z} = \bm{x}+\bm{y}$ is computed as
\begin{align}\label{eq:vector_sum}
\bm{z} = \sum_{i=1}^{N_{\mathrm{par}}}\bm{P}^{(i)} \bm{z}^{(i)} = \sum_{i=1}^{N_{\mathrm{par}}}\bm{P}^{(i)} \left\{\bm{x}^{(i)}+\bm{y}^{(i)}\right\}.
\end{align}
The parallel computation of the vector sum involves operations on the vectors $\bm{x}^{(i)}$ and $\bm{y}^{(i)}$ corresponding to the internal node set $V^{(i)}$ of subdomain $i$, and the vector sum ${\bm{z}}^{(i)}$ is computed in each subdomain. This operation can be performed without MPI communication.

For the vector inner product $\alpha=\bm{x}^{\mathrm{T}}\bm{y}$, the following relation holds.
\begin{equation}\label{eq:vecvec}
 \alpha= \sum_{i=1}^{N_{\mathrm{par}}}\alpha^{(i)}= \sum_{i=1}^{N_{\mathrm{par}}}{\bm{x}^{(i)}}^{\mathrm{T}}\bm{y}^{(i)}.
\end{equation}
Here, $\alpha^{(i)}$ is the value of the vector inner product associated with subdomain $i$. The parallel computation of the vector inner product is performed by first computing local inner products using the vectors $\bm{x}^{(i)}$ and $\bm{y}^{(i)}$ corresponding to the internal node set $V^{(i)}$ of subdomain $i$, and then summing the results over all parallel processes using a global communication function, namely the allreduce function.

Next, we describe the parallel computation of the matrix-vector product $\bm{z}=\bm{A}\bm{x}$. The parallel computation of the matrix-vector product is divided into the computation $\bm{A}^{(i,i)} \bm{x}^{(i)}$ associated with subdomain $i$ and the computation $\bm{A}^{(i,j)} \bm{x}^{(j)}$ associated with subdomain $i$ and its adjacent subdomains $j$, as expressed in Eq.~(\ref{eq:matvec}).
\begin{align}\label{eq:matvec}
\bm{z}& = \sum_{i=1}^{N_{\mathrm{par}}}\bm{P}^{(i)}\bm{z}^{(i)}=\sum_{i=1}^{N_{\mathrm{par}}} \bm{P}^{(i)}  \left\{ \bm{A}^{(i, i)} \bm{x}^{(i)}+\sum_{{j}\, \in\, \mathrm{nbhd}(\Tilde{v}_i, \Tilde{G})}
\bm{A}^{(i, j)} \bm{x}^{(j)}\right\}.
\end{align}
Here, $\mathrm{nbhd}(\Tilde{v}_i,\Tilde{G})$ is a function that returns all indices of nodes adjacent to the node $\Tilde{v}_i \in \Tilde{V}$ in the graph $\Tilde{G} = (\Tilde{V},\Tilde{E})$, and is defined as
\begin{equation}
\operatorname{nbhd}(\tilde{v}_i, \tilde{G}) := \left\{\, j \mid (\tilde{v}_j, \tilde{v}_i) \in \tilde{E},\ j \neq i \,\right\}.
\end{equation}

In the operations in Eqs.~(\ref{eq:vector_sum}) and (\ref{eq:matvec}), the computational results are stored only in the vector $\bm{z}^{(i)}$ associated with the internal node set $V^{(i)}$ of subdomain $i$. On the other hand, to compute $\bm{A}^{(i,j)}\bm{x}^{(j)}$, the vector $\bm{x}^{(j)}$ of the adjacent subdomain $j$ is required, but it is not stored by the process assigned to subdomain $i$. Therefore, the computation is performed by obtaining $\bm{x}^{(j)}$ through MPI communication from adjacent subdomains $j$ satisfying $j \in \mathrm{nbhd}(\tilde{v}_i,\tilde{G})$. The communication table required to obtain the vector $\bm{x}^{(j)}$ can be constructed based on the adjacency information of the metagraph $\Tilde{G}$. In practice, as illustrated in Fig.~\ref{fig:matvec}, to minimize the amount of data communication, only the vector components of $\bm{x}^{(j)}$ corresponding to the nonzero entries of the matrix $\bm{A}^{(i,j)}$ are communicated, based on the information of the subgraphs $G^{(i)}$ and $G^{(j)}$. The matrix $\bm{A}^{(i,i)}$ corresponds to the internal node set $V^{(i)}$ of subdomain $i$ and can be generated independently in each subdomain. Similarly, the matrix $\bm{A}^{(i,j)}$ can also be generated independently in each subdomain without MPI communication because the overlap region $V^{(i)}_{\mathrm{ovl}}$ of subdomain $i$ contains the graph information required to compute $\bm{A}^{(i,j)}$. Thus, the matrix-vector product to be computed by each parallel process is defined based on the subgraph $G^{(i)}$ and the metagraph $\Tilde{G}$.

\subsection{Hierarchical Metagraph-Based Parallelization}

\begin{figure*}[t]
    \begin{center}
    \includegraphics[scale=0.46]{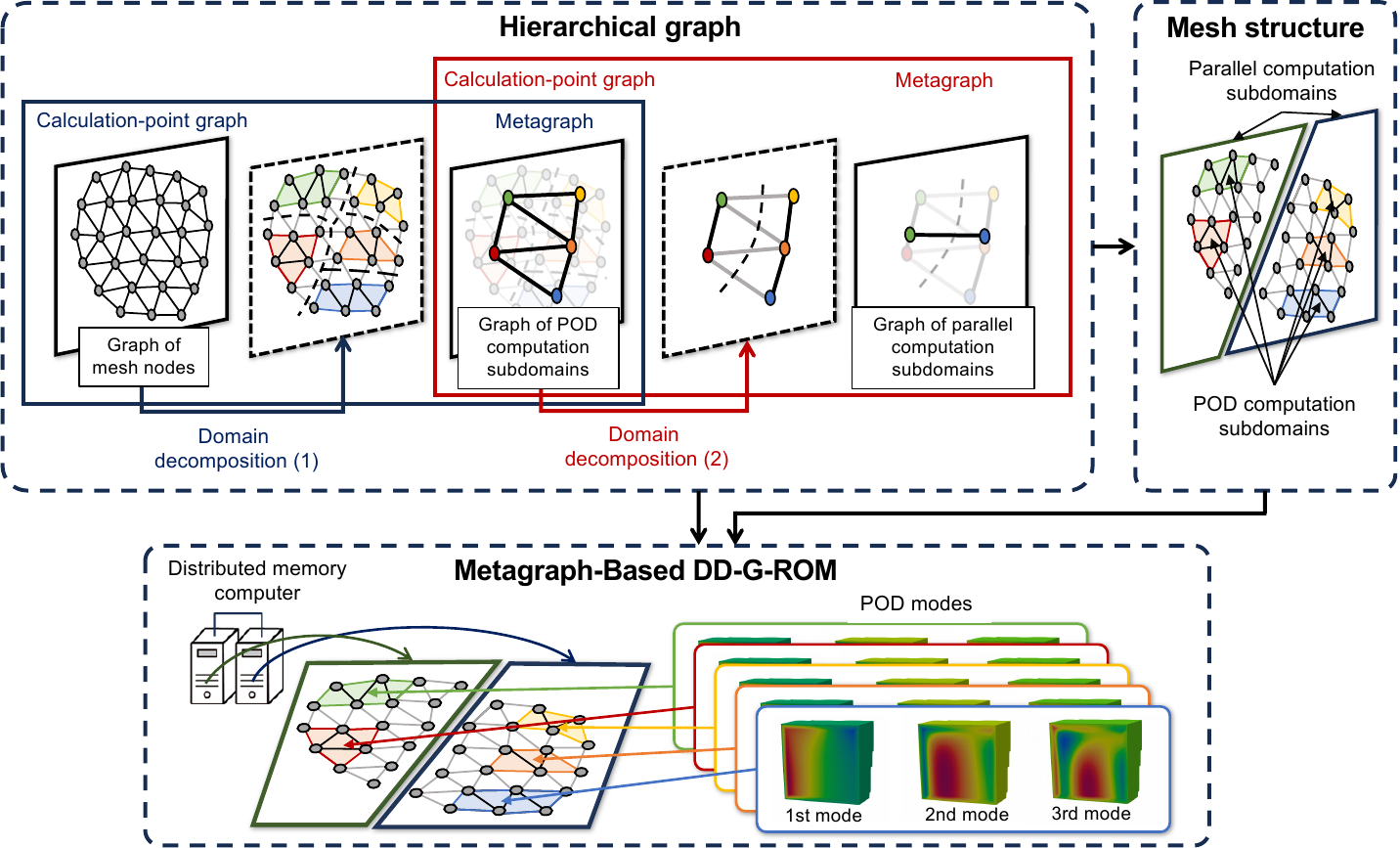}
    \end{center}
    \vspace{-10pt}
    \caption{\baselineskip=13pt Schematic illustration of the metagraph-based DD-G-ROM framework: hierarchical graph representation of POD computation subdomains and parallel computation subdomains.}
    \label{fig:main}
\end{figure*}

In this study, by exploiting graph structures, L-POD is generalized so that it can be computed on arbitrary subdomains. In addition, we propose a hierarchical domain-decomposition method based on a metagraph structure to independently define the ``subdomains used to compute POD bases'' in L-POD, introduced in Section~\ref{syou:Local POD}, and the ``subdomains used for domain-decomposition-based parallelization,'' introduced in Section~\ref{syou:parallel_comp}, thereby realizing distributed-memory parallel computation of L-POD. For conciseness, hereafter we refer to the subdomains used to compute POD bases as POD computation subdomains and the subdomains used for domain-decomposition-based parallelization as parallel computation subdomains. This study focuses on problems in which the number of POD computation subdomains is larger than the number of parallel computation subdomains. 
Figure~\ref{fig:main} schematically illustrates the proposed MBDD-G-ROM and its graph hierarchy of POD computation subdomains and parallel computation subdomains.
The proposed two-level domain-partitioning strategy consists of two successive partitioning stages. In the first partitioning stage, the finite element graph is partitioned to define the POD computation subdomains. In the second partitioning stage, the metagraph constructed from these subdomains is partitioned to define the parallel computation subdomains.

In the first partitioning stage, based on the definitions of graphs and domain decomposition introduced in Section~\ref{syou:DD}, the graph corresponding to the finite element nodes is partitioned (Domain decomposition 1 in Fig.~\ref{fig:main}) to generate POD computation subdomains. The POD computation subdomains generated in this manner can be interpreted as coarse finite elements to which low-dimensional basis functions are assigned by L-POD. Usually, the basis functions generated by POD are modes defined over the entire domain, and their graph structure is all-to-all connected, corresponding to a dense matrix. Consequently, in the absence of a sufficiently sparse block structure, domain-decomposition-based data distribution leads to extensive interprocess communication, making it highly challenging to achieve favorable strong-scaling performance on distributed-memory systems. In contrast, L-POD divides the domain and has a sparse block matrix structure corresponding to the metagraph described above, which enables model reduction and domain-decomposition-based parallel analysis to be combined.

In the second partitioning stage, according to the definition of the metagraph introduced in Section~\ref{syou:DD}, the metagraph $\tilde{G}$ corresponding to the set of POD computation subdomains is constructed. At this stage, if there exists a point $\bm{x}\ \in \mathbb{R}^3$ at which the product of the POD basis $\bm \Phi^{(a)}(\bm{x})$ of metanode $a$ and the POD basis of metanode $b$, $\bm \Phi^{(a)}(\bm{x})^{\mathrm{T}} \bm \Phi^{(b)}(\bm{x})$, is nonzero, namely if
\begin{equation}\label{eq:def_edge_meta}
\exists\, \bm{x} \in \mathbb{R}^3, \ \bm \Phi^{(a)}(\bm{x})^{\mathrm{T}} \bm \Phi^{(b)}(\bm{x}) \neq \bm 0
\end{equation}
holds, metanode $a$ is regarded as adjacent to metanode $b$. By applying graph partitioning again to the resulting metagraph (Domain decomposition 2 in Fig.~\ref{fig:main}), the parallel computation subdomains can be defined. As illustrated in the lower part of Fig.~\ref{fig:main}, each parallel computation subdomain contains a specific set of POD computation subdomains, and distributed-memory parallel computation is realized by applying the domain-decomposition-based parallel computation described in Section~\ref{syou:parallel_comp}. This second-stage graph partitioning is performed to define the parallel computation subdomains and does not repartition the finite element node set $V$ or the POD bases $\bm{\Phi}^{(i)}$. Therefore, the numerical values of $\bm \Phi^{(i)}$, $\tilde{\bm J}^{(i,j)}$, and $\tilde{\bm r}^{(i)}$ remain unchanged.

The following describes the specific matrix structure obtained when the proposed method is applied to L-POD. First, consider the metagraph $\tilde{G}=(\tilde{V},\tilde{E})$ obtained by partitioning the finite element graph (upper part of Fig.~\ref{fig:main}). The global matrix of L-POD corresponding to the metagraph $\tilde{G}$ is given by
\begin{equation}\label{eq:LocalPOD1}
  \tilde{\bm{J}}=
  \begin{bmatrix}
    \tilde{\bm{J}}_{1,1} & \tilde{\bm{J}}_{1,2} & \cdots & \tilde{\bm{J}}_{1,|\tilde V|}\\
    \tilde{\bm{J}}_{2,1} & \tilde{\bm{J}}_{2,2} & \cdots & \tilde{\bm{J}}_{2,|\tilde V|}\\
    \vdots & \vdots & \ddots & \vdots\\
    \tilde{\bm{J}}_{|\tilde V|,1} & \tilde{\bm{J}}_{|\tilde V|,2} & \cdots & \tilde{\bm{J}}_{|\tilde V|,|\tilde V|}
  \end{bmatrix}.
\end{equation}
Here, each block $\Tilde{\bm{J}}_{i,j}\in\mathbb{R}^{n_{\mathrm{POD}}^{(i)}\times n_{\mathrm{POD}}^{(j)}}$ represents the action from metanode $\tilde{v}_j$ to metanode $\tilde{v}_i$, and $\Tilde{\bm{J}}_{i,j}$ is the zero matrix when $(\tilde{v}_i,\tilde{v}_j)\notin \tilde{E}$.

Next, the metagraph $\tilde{G}$ is graph-partitioned (Domain decomposition 2 in Fig.~\ref{fig:main}), and the action of the internal node set $\tilde{V}^{(j)}$ of parallel computation subdomain $j$ on the internal node set $\tilde{V}^{(i)}$ of parallel computation subdomain $i$ is expressed as the following block matrix $\tilde{\bm{J}}^{(i,j)}$.
\begin{equation}
\tilde{\bm{J}}^{(i,j)}=
\begin{bmatrix}
\tilde{\bm{J}}^{(i,j)}_{1,1} & \tilde{\bm{J}}^{(i,j)}_{1,2} & \cdots & \tilde{\bm{J}}^{(i,j)}_{1,|\tilde V^{(j)}|}\\
\tilde{\bm{J}}^{(i,j)}_{2,1} & \tilde{\bm{J}}^{(i,j)}_{2,2} & \cdots & \tilde{\bm{J}}^{(i,j)}_{2,|\tilde V^{(j)}|}\\
\vdots & \vdots & \ddots & \vdots\\
\tilde{\bm{J}}^{(i,j)}_{|\tilde V^{(i)}|,1} & \tilde{\bm{J}}^{(i,j)}_{|\tilde V^{(i)}|,2} & \cdots & \tilde{\bm{J}}^{(i,j)}_{|\tilde V^{(i)}|,|\tilde V^{(j)}|}
\end{bmatrix}.
\end{equation}
Here, $\tilde{\bm{J}}^{(i,j)}_{p,q} \in
\mathbb{R}^{\,n_{\mathrm{POD}}^{(\operatorname{idx}(i,p))}\times
n_{\mathrm{POD}}^{(\operatorname{idx}(j,q))}}$ represents the action from the $q$th metanode $\tilde{v}_{\mathrm{idx}(j,q)}\in \tilde{V^{(j)}}$ of $\tilde V^{(j)}$ to the $p$th metanode $\tilde{v}_{\mathrm{idx}(i,p)}\in \tilde{V^{(i)}}$ of $\tilde V^{(i)}$, and $\tilde{\bm{J}}^{(i,j)}_{p,q}$ is the zero matrix when $(\tilde{v}_{\mathrm{idx}(j,q)},\tilde{v}_{\mathrm{idx}(i,p)})\notin \tilde{E}$.

Finally, using the matrix $\bm{P}^{(i)}\in
\mathbb{R}^{\sum_{k\in\tilde{\mathcal V}} n^{(k)}_{\mathrm{POD}}\times
\sum_{l\in\tilde{\mathcal V}^{(i)}} n^{(l)}_{\mathrm{POD}}}$, the relation between the matrices $\tilde{\bm{J}}$ before partitioning and $\tilde{\bm{J}}^{(i,j)}$ after partitioning is expressed as
\begin{equation}
\tilde{\bm{J}}=\sum_{i=1}^{N_{\mathrm{par}}}\sum_{j=1}^{N_{\mathrm{par}}}\bm{P}^{(i)}\,\tilde{\bm{J}}^{(i,j)}\,{\bm{P}^{(j)}}^{\mathrm{T}}.
\end{equation}
Using the $\mathrm{idx}$ function defined in Section~\ref{syou:DD}, the components of the matrix $\bm P$ are given by
\begin{equation}
\bm{P}^{(i)}_{m,\,n}=
\begin{cases}
\bm{I} \in \mathbb{R}^{\,n_{\mathrm{POD}}^{(m)}\times n_{\mathrm{POD}}^{(m)}}, & m=\operatorname{idx}(i,n)\\
\bm{0} \in \mathbb{R}^{\,n_{\mathrm{POD}}^{(m)}\times n_{\mathrm{POD}}^{\left(\operatorname{idx}(i,n)\right)}}, & \text{otherwise}
\end{cases}.
\end{equation}

\subsection{Metanode Weight Selection for Parallel Computation}\label{syou:staticlb}

Because the number of basis vectors stored in each POD computation subdomain differs in this study, the computational load can become imbalanced among parallel processes, leading to a reduction in the overall computational efficiency. To address this issue, we introduce load balancing based on a weighted graph.

First, a metagraph $\tilde{G} = (\tilde{V}, \tilde{E})$ is constructed by treating the POD computation subdomains as nodes and representing the adjacency relations among POD computation subdomains as edges. At this stage, each metanode $\tilde{v}_i$ is assigned a node weight $n_{\mathrm{weight}}^{(i)}$, which represents the estimated computational cost of the corresponding subdomain. If necessary, the estimated communication volume between adjacent POD computation subdomains can be assigned as an edge weight; however, in this study, all edge weights are set to one.

Next, this weighted graph is partitioned into subgraphs $\tilde{G}^{(i)}\ (i=1,2,\ldots,N_{\mathrm{par}})$. The graph partitioning is performed so as to balance the sum of the node weights contained in each internal node set $\tilde{V}^{(i)}$, while minimizing the total weight of the edges crossing the internal node sets $\tilde{V}^{(i)}$.

The proposed method described above does not depend on how the POD computation subdomains are partitioned and is applicable to arbitrary domain decompositions. Various methods have been used in previous studies for domain decomposition, including regular decomposition based on structured grids \cite{Pekurovsky2012P3DFFT} and decomposition using space-filling curves \cite{Bader2013SFC}. However, the focus of this study is not to propose a domain-decomposition method itself, but to assign the partitioned POD computation subdomains to parallel computation subdomains while considering the computational load based on the number of basis vectors. Therefore, in this study, we adopt the same graph-partitioning-based method for partitioning both the POD computation subdomains and the parallel computation subdomains. Specifically, the partitioning is determined so as to balance the computational load by equalizing the sum of the node weights contained in each subdomain, while minimizing the sum of the edge weights crossing between subdomains, namely the edge cut. To solve this weighted graph partitioning problem, this study uses the graph partitioning library METIS \cite{Karypis1998METIS}.

%% file: include/numerical_analysis.tex
\section{Numerical Examples}
This section first considers an unsteady diffusion equation, which is a linear partial differential equation, and quantitatively evaluates the accuracy, parallel efficiency, and computational efficiency of the proposed method in Section~5.1. An example of static load balancing is also presented. Section~5.2 then presents an analysis of the Navier--Stokes equations using MBDD-G-ROM to verify the applicability of the proposed method to nonlinear partial differential equations.

\begin{table}[t]
\small
\centering
\caption{Specifications of the SQUID supercomputer.}
\begin{tabular}[t]{ll}
\hline
Component&\ \;Specification\\
\hline
\hline
CPU&\begin{tabular}{l}Intel Xeon Platinum 8368\\
(Ice Lake, 2.4 GHz, 38 cores) $\times 2$\end{tabular}
\\
Memory&\ \;256 GB\\
Peak performance&\ \;5.84 TFLOPS\\
Number of nodes&\ \;1,520\\
\hline
\end{tabular}\label{tablesquid}
\vspace{-12pt}
\end{table}

\subsection{Unsteady Diffusion Equation}

The distributed-memory computer used in this verification was the SQUID supercomputer system at Osaka University \cite{squid}. The system specifications are listed in Table~\ref{tablesquid}.

\begin{figure*}[t]
    \vspace{10pt}
    \begin{center}
    \includegraphics[scale=0.6]{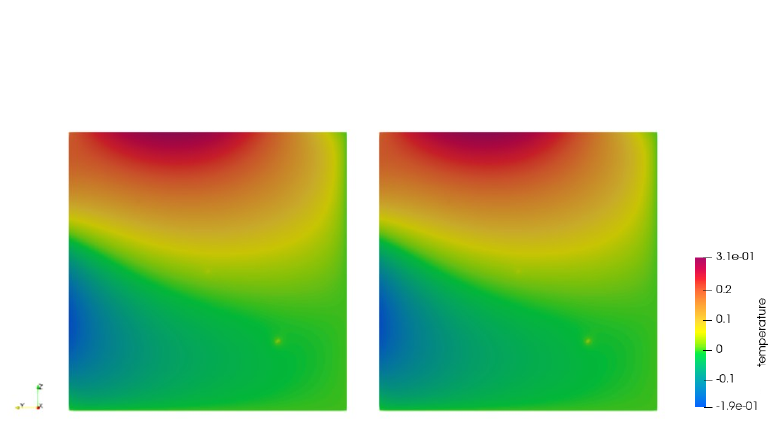}
    \end{center}
    \vspace{-5pt}
    \caption{\baselineskip=13pt Contours of the numerical solution at $t=10s$ on a cross-section normal to the $x$-axis; FOM (left) and ROM (right).}
    \vspace{-5pt}
    \label{fig:compare_FEM_L-POD}
\end{figure*}

In this verification, the following unsteady diffusion equation is used as the governing equation:
\begin{equation}\label{eq:poisson_equ}
    \frac{\partial u}{\partial t}-\,\nabla\cdot ( k\,\nabla u)=Q.
\end{equation}
Here, $k$ is the diffusion coefficient, $u$ is the unknown variable representing the diffusing physical quantity, and $Q$ is the source term representing the physical quantity generated in the analysis domain $\Omega$, where $\Omega$ is a bounded three-dimensional domain. The boundary conditions are given by
\begin{align}
    u&=u_D\quad \mathrm{on} \quad \mathrm{\Gamma}_g,\label{eq:dirichlet_bc}\\
    \nabla {u} \cdot \bm{n}&= {q}_N \quad \mathrm{on} \quad \mathrm{\Gamma}_h,\label{eq:neumann_bc}
\end{align}
where $u_D$ is the value of the solution $u$ on the Dirichlet boundary $\mathrm{\Gamma}_g$, $\bm{n}$ is the outward normal vector on the Neumann boundary $\mathrm{\Gamma}_h$, and $q_N$ is the prescribed function representing the physical quantity on the boundary $\mathrm{\Gamma}_h$. In this verification, the diffusion coefficient is set to $k(x, y, z)=1+x^3+y^3+z^3$. The source term is defined as the sum of point heat sources, and its time dependence is given by
\begin{equation}
q(t)=1000\bigl(1+\sin t\bigr).
\end{equation}
Using the set of point-source locations
\begin{equation}
\mathcal{S}=\{(2.5,3.75,3.75),\ (2.5,2.75,2.5),\ (2.5,1.25,1.25)\},
\end{equation}
the source term of the unsteady diffusion equation can be expressed as
\begin{equation}
Q(\boldsymbol{x},t)
= q(t)\sum_{\boldsymbol{x}_s\in\mathcal{S}}\delta(\boldsymbol{x}-\boldsymbol{x}_s),
\qquad
\delta(\boldsymbol{x}-\boldsymbol{x}_s)
=\delta(x-x_s)\,\delta(y-y_s)\,\delta(z-z_s),
\end{equation}
where \(\delta\) denotes the Dirac delta distribution. In the implementation, this term was applied as point heat sources at mesh nodes. The following Dirichlet condition was imposed on the entire boundary $\mathrm{\Gamma}$:
\begin{align}\label{eq:settings_diri}
    u &= \sin\left(0.25x\right)\sin\left(0.5y\right) \sin z \sin t \quad \mathrm{on} \quad \mathrm{\mathrm{\Gamma}}.
\end{align}

In this verification, common spatial and temporal discretization settings were used for the numerical method employed to generate the basis in the offline phase and for the numerical method to be reduced in the online phase. The finite element method was used for the spatial discretization, and the implicit Euler method was used for the temporal discretization. The analysis domain $\Omega$ was set to the three-dimensional cubic domain $[0,5]^3$, and the mesh was constructed so that all elements were cubic. First-order Lagrange interpolation functions were used as the finite element basis functions, and the time-step size $\Delta t$ was set to $0.01$s. In the online phase, ROMs were constructed and analyzed using the POD-Galerkin method.

For both the finite element and ROM analyses, MONOLIS (Metagraph-Oriented Network for Linear Iterative Solvers) \cite{monolis} was used for matrix computations. The CG (conjugate gradient) method was used to solve the systems of linear equations, with a convergence tolerance of $1.0\times 10^{-9}$ and a maximum of 10,000 iterations. Diagonal scaling was used as the preconditioner. The performance of the diagonal scaling preconditioner is independent of the number of subdomains. 
ScaLAPACK \cite{ScaLAPACK} was used to perform the distributed-memory parallel singular value decomposition.

\subsubsection{Verification of the Computational Accuracy of L-POD}
This subsection verifies the computational accuracy when the number of basis vectors and the number of POD computation subdomains are varied, with the aim of evaluating the accuracy in serial computations. 

First, finite element analyses were performed over the initial 100 time steps, and the snapshot matrix was constructed from the solution vectors obtained during this process. The snapshot matrix was then subjected to singular value decomposition to generate the POD basis used in the ROM analysis. The number of POD basis vectors used in each POD computation subdomain was determined according to Eq.~(\ref{eq:base_sel}). The same computational settings were also used in Section~\ref{kensyo:calceff}. Subsequently, ROM analyses were performed using the generated basis over the range from time step 101 to time step 1,000, namely from $t=1.01\ \mathrm{s}$ to $t=10.0\ \mathrm{s}$.

For error evaluation, the numerical solution $u_{\mathrm{ROM}}$ obtained by ROM analysis was compared with the numerical solution $u_{\mathrm{FEM}}$ obtained by finite element analysis using the relative error based on the $L^2$ norm, defined in Eq.~(\ref{eq:settings_error}).
\begin{equation}\label{eq:settings_error}
\varepsilon_{\mathrm{ROM}}
=\frac{\sqrt{\int_{\Omega}|u_{\mathrm{ROM}}-u_{\mathrm{FEM}}|^2d \Omega}}{\sqrt{\int_{\Omega}|u_{\mathrm{FEM}}|^2d \Omega}}.
\end{equation}

\begin{figure}[t]
    \begin{center}
    \includegraphics[scale=0.44]{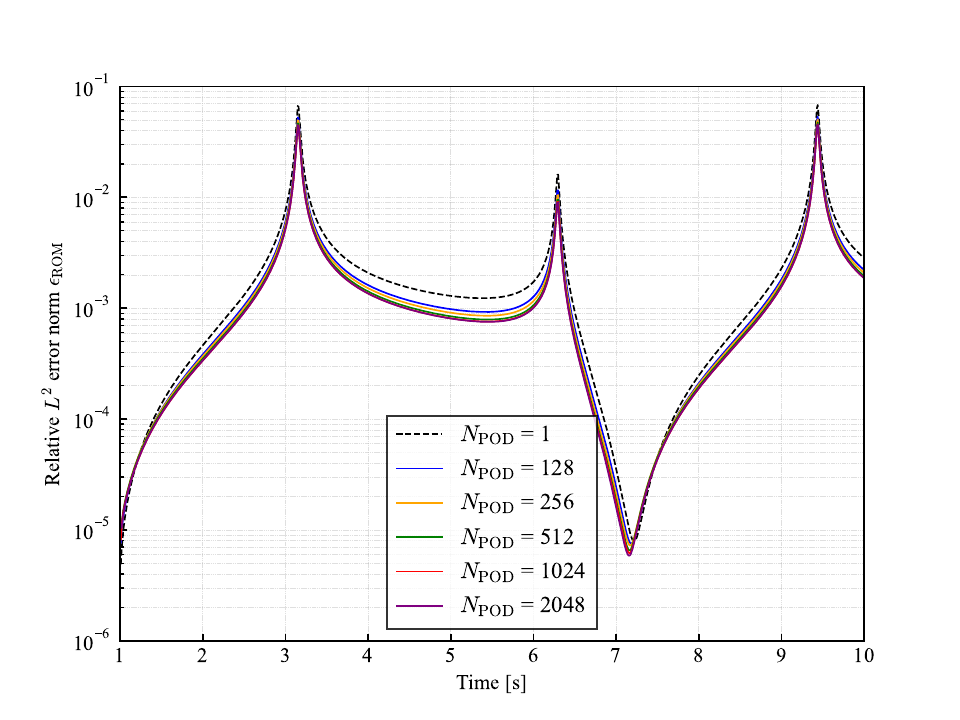}
    \end{center}
    \vspace{-10pt}
    \caption{\baselineskip=13pt Relative $L^2$ error norm of L-POD for different numbers of subdomains.}
    \vspace{-5pt}
    \label{fig:L2error_num_subdomains}
\end{figure}

Six cases with 1, 128, 256, 512, 1,024, and 2,048 POD computation subdomains were examined for a finite element analysis with 1,030,301 degrees of freedom. The threshold for basis selection was set to $\varepsilon_{\mathrm{POD}}=10^{-6}$. The total number of POD basis vectors in the POD computation subdomains is shown in Table~\ref{tab:values_num_modes}. Figure~\ref{fig:compare_FEM_L-POD} shows contours of the numerical solutions at 10s on a cross-section normal to the $x$-axis ($x$=2.5) for the finite element and ROM analyses, and Fig.~\ref{fig:L2error_num_subdomains} shows the time history of the relative error in Eq.~(\ref{eq:settings_error}). Table~\ref{tab:values_num_modes} and Fig.~\ref{fig:L2error_num_subdomains} show that increasing the number of subdomains increases the number of basis functions and improves the computational accuracy. These results demonstrate the validity of the implementation.

\begin{table}[tbp]
\centering
\small
\caption{Total number of POD basis vectors for different numbers of POD computation subdomains.}
\label{tab:values_num_modes}
\begin{tabular}{lcccccc}
\hline
Number of POD computation subdomains & 1 & 128 & 256 & 512 & 1,024 & 2,048 \\
\hline
Total number of POD basis vectors & 8 & 726 & 1,331 & 2,455 & 4,619 & 8,702 \\
\hline
\end{tabular}
\end{table}

\begin{figure*}[t]
  \centering
  \begin{subfigure}[t]{0.4\linewidth}
    \hspace{25pt}
    \includegraphics[scale=0.5]{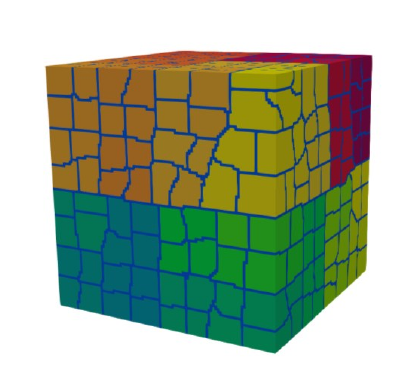}
    \subcaption{\small POD computation subdomains.}
  \end{subfigure}
  \hspace{0.05\linewidth}
  \begin{subfigure}[t]{0.4\linewidth}
    \hspace{25pt}
    \includegraphics[scale=0.5]{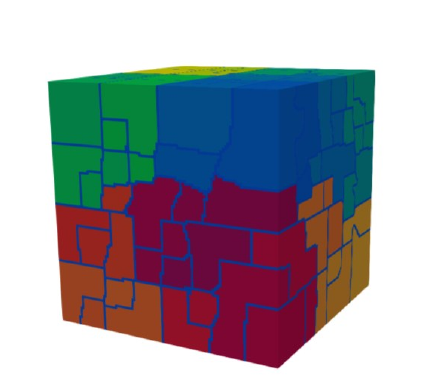}

    \subcaption{\small Parallel computation subdomains.}
  \end{subfigure}
  \vspace{5.5pt} 
    \caption{\baselineskip=13pt Example of the two-level decomposition into POD computation and parallel computation subdomains (degrees of freedom = 1,030,301, $N_{\mathrm{POD}}=512$, $N_{\mathrm{par}}=128$).}
    \vspace{-5pt}
    \label{fig:DDM_mesh}
\end{figure*}

\subsubsection{Verification of the Computational Efficiency of the Proposed Method}\label{kensyo:calceff}
This subsection evaluates the parallel computational efficiency of the proposed method through strong-scaling tests, in which the number of parallel computation subdomains is varied while the number of POD computation subdomains is fixed. As discussed in Section~\ref{syou:staticlb}, node weights can be assigned to the metagraph for static load balancing. In this verification, however, all node weights are set to unity to evaluate the baseline parallel efficiency of the proposed method without load balancing.

\begin{figure*}[t]
  \centering
  \begin{subfigure}[t]{0.4\linewidth}
    \includegraphics[scale=0.44]{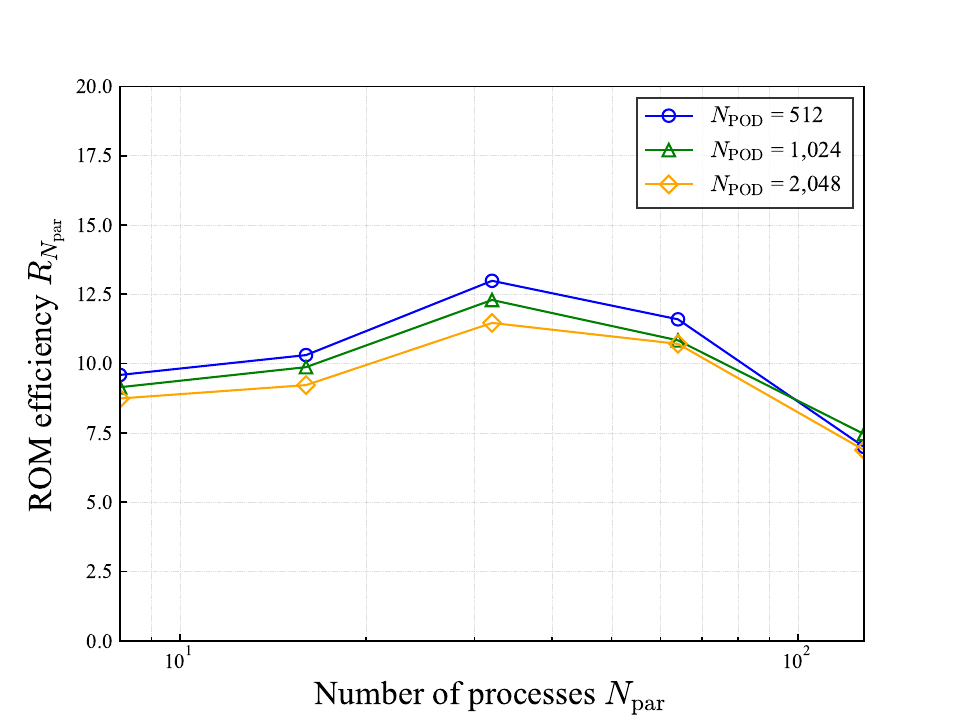}
    \subcaption{Degrees of freedom = 1,030,301.}
  \end{subfigure}
  \hspace{0.05\linewidth}
  \begin{subfigure}[t]{0.4\linewidth}
    \includegraphics[scale=0.44]{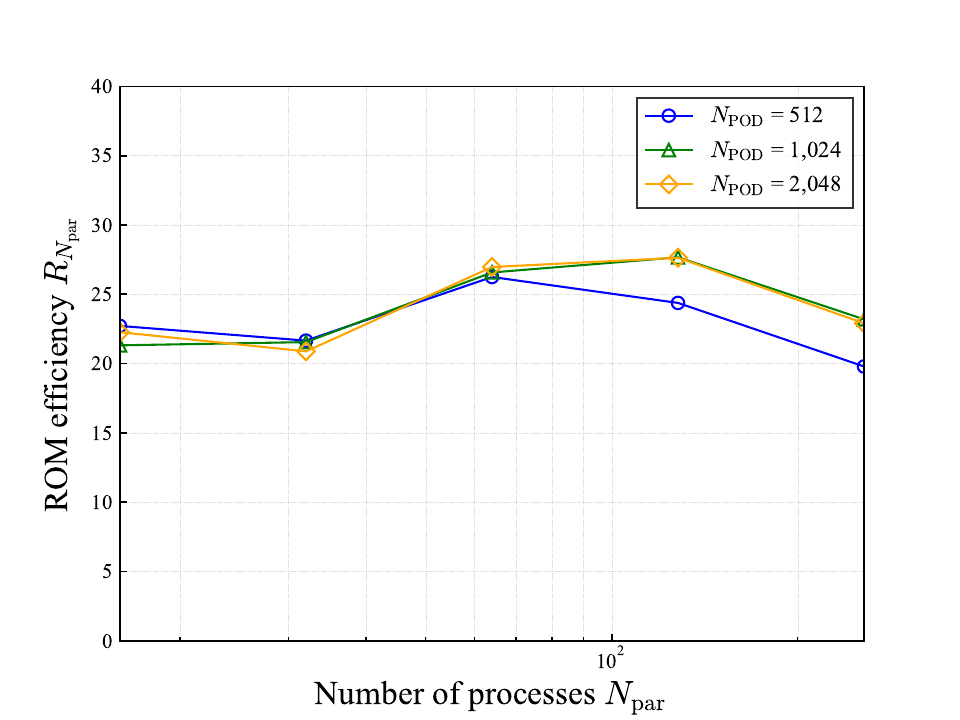}
    \subcaption{Degrees of freedom = 10,077,696.}
  \end{subfigure}
  \vspace{5.5pt} 
  \caption{ROM efficiency based on total runtime for different numbers of parallel processes.}
\label{verifi:strong_ROM_efficiency}
\end{figure*}

Three cases with 512, 1,024, and 2,048 POD computation subdomains were examined for two finite element models with 1,030,301 and 10,077,696 degrees of freedom. For the case with 1,030,301 degrees of freedom, the number of parallel processes was varied as 1, 8, 16, 32, 64, 128, and 256. For the case with 10,077,696 degrees of freedom, the number of parallel processes was varied as 16, 32, 64, 128, and 256. The resulting domain decomposition is shown in Fig.~\ref{fig:DDM_mesh}. To keep the effect of memory bandwidth limitation approximately constant, the maximum number of processes per node was set to 16. Thus, the number of compute nodes used was set to 1 for 1, 8 and 16 parallel processes, and to 2, 4, 8, and 16 for 32, 64, 128, and 256 parallel processes, respectively. 
In the basis-selection criterion given by Eq.~\eqref{eq:base_sel}, the same threshold, $\varepsilon_{\mathrm{POD}}=10^{-6}$, was used for all POD computation subdomains.
Each measurement was performed three times for each case, and the average value was used for evaluation. In the case with 10,077,696 degrees of freedom, the verification problem could not be executed serially because of memory capacity limitations. Therefore, the serial computation time $T_1$ was estimated by applying Amdahl's law based on the results obtained with 32 and 64 parallel processes using two and four compute nodes, respectively. This computational example demonstrates that the proposed distributed-memory parallel computing method alleviates memory capacity limitations and confirms the effectiveness of the proposed method.

To evaluate the parallel computing performance, the ROM efficiency $R_{N_{\mathrm{par}}}$ and the speed-up factor $S_{N_{\mathrm{par}}}$ were used as performance metrics. The ROM efficiency represents the acceleration factor of the ROM analysis relative to the finite element analysis, and the two metrics are defined as
\begin{equation}\label{eq:ROM_efficiency}
R_{N_{\mathrm{par}}}=\frac{T^{\mathrm{FEM}}_{N_{\mathrm{par}}}}{T^{\mathrm{ROM}}_{N_{\mathrm{par}}}},
\end{equation}
\begin{equation}\label{eq:speed-upfactor}
    S_{N_{\mathrm{par}}}=\frac{T_1}{T_{N_{\mathrm{par}}}}.
\end{equation}
Here, $T_{N_{\mathrm{par}}}^{\mathrm{FEM}}$ and $T_{N_{\mathrm{par}}}^{\mathrm{ROM}}$ denote the computation times using the finite element method and ROM, respectively, with $N_{\mathrm{par}}$ parallel processes, and $T_{N_{\mathrm{par}}}$ denotes the computation time with $N_{\mathrm{par}}$ parallel processes.

\begin{table*}[tbp]
\small
    \centering
   \caption{Total number of POD basis vectors for different numbers of POD computation subdomains and degrees of freedom.}
    \begin{tabular}{cc|ccccc}\hline
    \multicolumn{2}{c|}{Number of POD computation subdomains} & 512 & 1,024 & 2,048 \\ \hline
    \multirow{2}{*}{\makecell{Degrees of freedom}}
& \ \:\:1,030,301 & 2,455 & 4,619 & 8,702&  \\
& 10,077,696 & 2,416 & 4,546 & 8,586 \\\hline
    \end{tabular}
    \label{tab:values_strong1}
\end{table*}

First, the total computation time was measured for 100 time steps, from time step 101 to time step 200, in each analysis case. Figure~\ref{verifi:strong_ROM_efficiency} shows the ROM efficiency evaluated using the average total computation time per time step. Figure~\ref{verifi:strong_ROM_efficiency} shows that the ROM efficiency is largely maintained with respect to the number of parallel processes in all analysis cases, demonstrating that the proposed method achieves good parallel computing performance. The ROM efficiency also tends to increase for cases with larger finite element degrees of freedom. This tendency can be interpreted as follows. The total number of POD basis vectors in the POD computation subdomains is listed in Table~\ref{tab:values_strong1}. As shown in Table~\ref{tab:values_strong1}, the number of reduced degrees of freedom in the ROM analysis is almost constant regardless of the number of degrees of freedom in the finite element analysis. Therefore, for cases with a larger number of finite element degrees of freedom, the reduction ratio in the number of degrees of freedom achieved by the ROM analysis becomes relatively larger. For these reasons, the ROM efficiency is considered to increase as the number of degrees of freedom in the finite element analysis becomes larger.

Next, we assess the strong-scaling performance of the proposed method in terms of speed-up.
For this assessment, the ROM runtime at each time step is divided into two components:
\begin{itemize}
\item[(a)] reduced-system assembly, including the construction of the reduced Jacobian and residual vector;
\item[(b)] solution of the reduced Newton system in Eq.~\eqref{eq:NR_reduced}.
\end{itemize}
The total ROM runtime is defined as the sum of these two components.
For the unsteady diffusion equation~\eqref{eq:poisson_equ}, which is linear, the reduced coefficient matrix
$\bm{\Phi}^{\mathrm{T}}\bm{J}\bm{\Phi}$ is assembled only once for each analysis case.
Because this one-time setup cost accounts for only a small fraction of the overall runtime, it is excluded from the strong-scaling assessment.
Consequently, component~(a) for the diffusion problem includes only the per-time-step assembly of the reduced right-hand-side vector.

\begin{figure*}[htbp]
  \centering
  \begin{subfigure}[t]{0.4\linewidth}
    \includegraphics[scale=0.44]{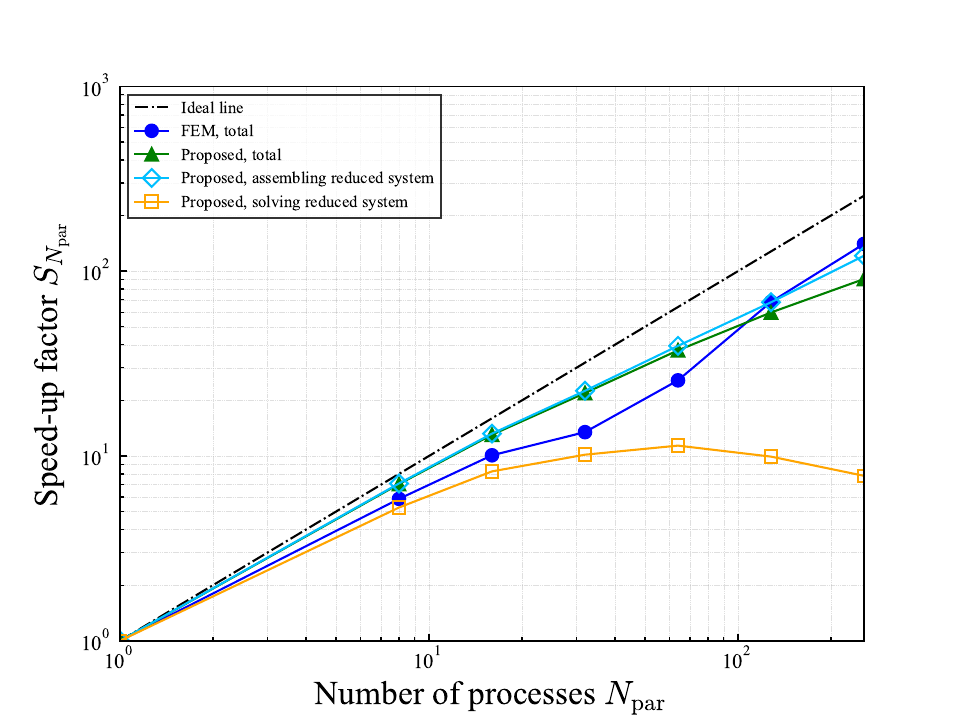}
    \vspace{-1em}
    \subcaption{$N_{\mathrm{POD}}= 512$.}
  \end{subfigure}
  \hspace{0.05\linewidth}
  \begin{subfigure}[t]{0.4\linewidth}
    \includegraphics[scale=0.44]{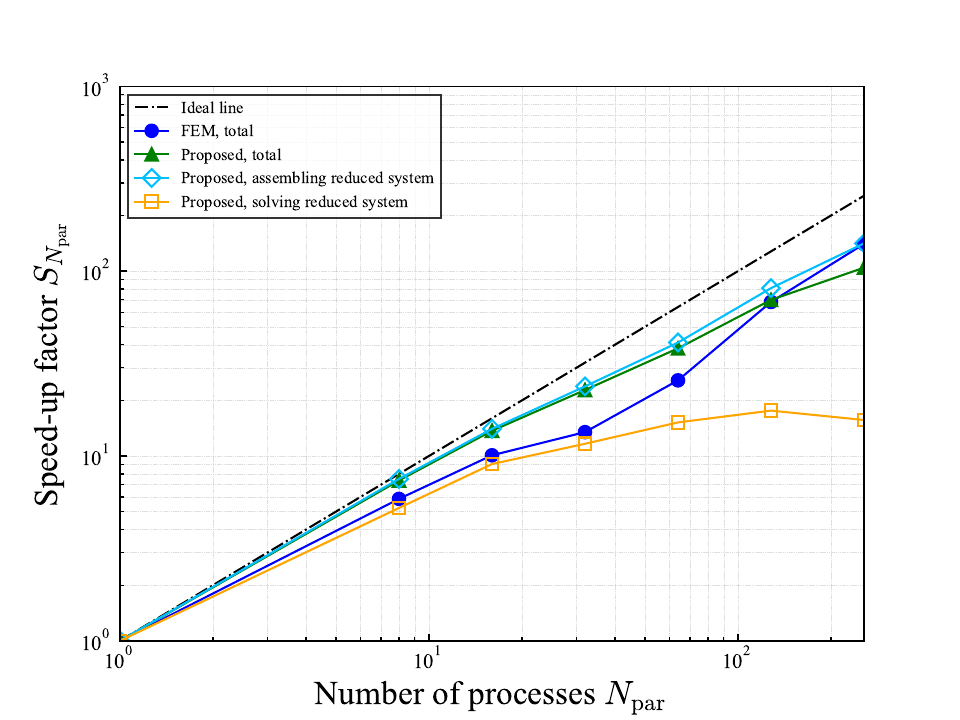}
        \vspace{-1em}
    \subcaption{$N_{\mathrm{POD}}$ = 1,024.}
  \end{subfigure}
  \vspace{0.01\linewidth}
          \vspace{0.5em}
  \begin{subfigure}[t]{0.4\linewidth}
    \includegraphics[scale=0.44]{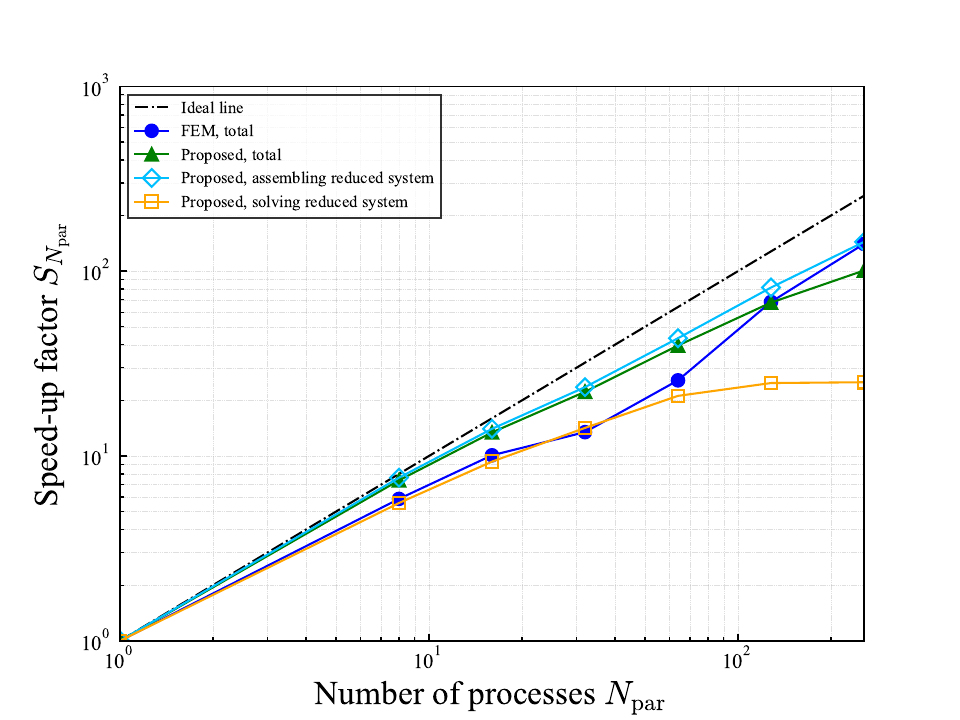}
        \vspace{0.5em} 
    \subcaption{$N_{\mathrm{POD}}$ = 2,048.}
  \end{subfigure}
  \caption{Speed-up factors of the conventional FEM and the proposed method (Degrees of freedom = 1,030,301).}
  \label{verifi:strong_speed-up_small}
\end{figure*}

\begin{figure*}[htbp]
  \centering
  \begin{subfigure}[t]{0.4\linewidth}
    \includegraphics[scale=0.44]{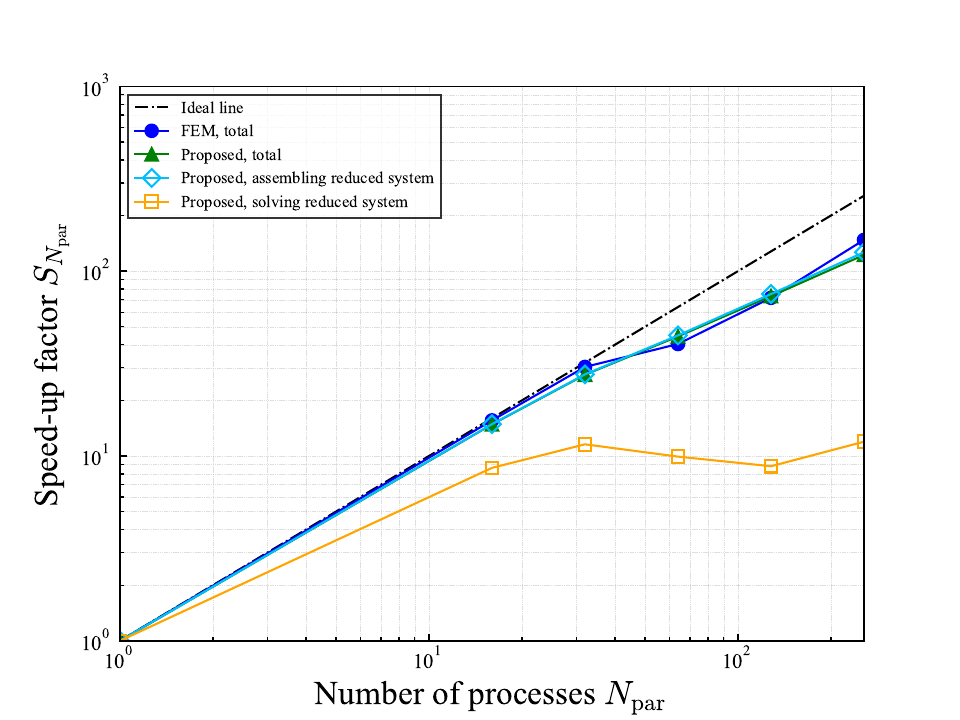}
    \vspace{-1em}
    \subcaption{$N_{\mathrm{POD}}= 512$.}
  \end{subfigure}
  \hspace{0.05\linewidth}
  \begin{subfigure}[t]{0.4\linewidth}
    \includegraphics[scale=0.44]{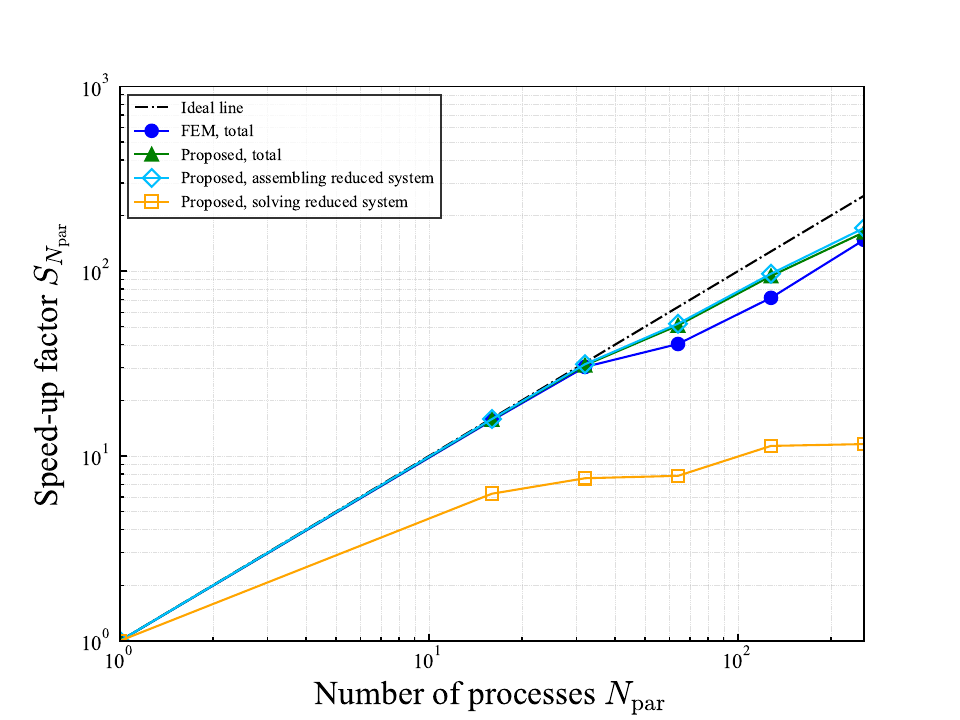}
    \vspace{-1em}
    \subcaption{$N_{\mathrm{POD}}$ = 1,024.}
  \end{subfigure}
  \vspace{0.01\linewidth}
            \vspace{0.5em}
  \begin{subfigure}[t]{0.4\linewidth}
    \includegraphics[scale=0.44]{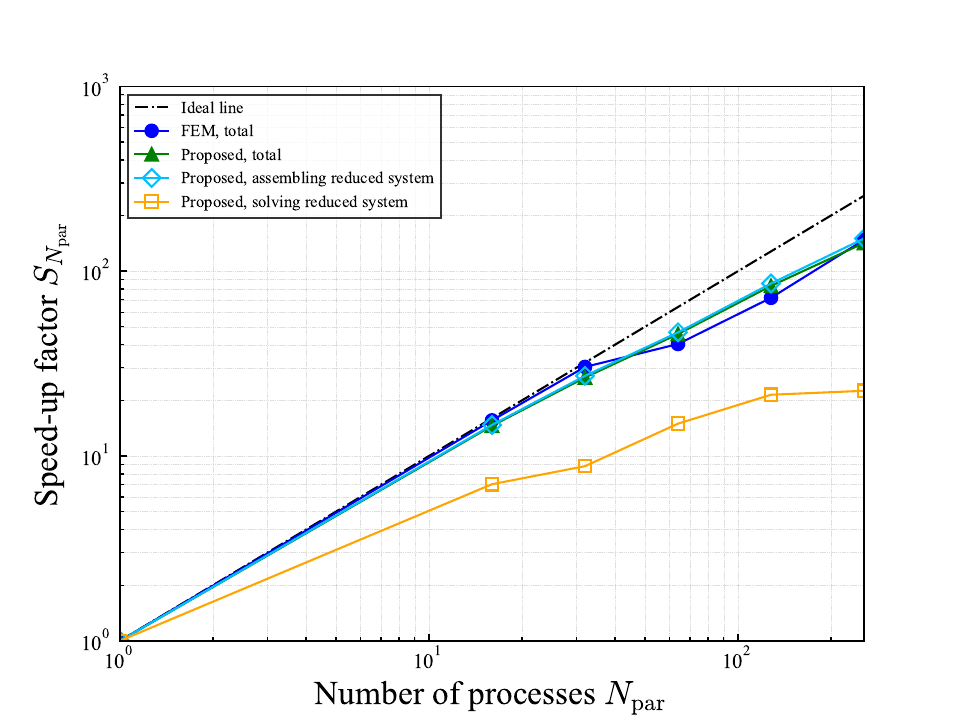}
\vspace{0.5em}
    \subcaption{$N_{\mathrm{POD}}$ = 2,048.}
  \end{subfigure}
  \caption{Speed-up factors of the conventional FEM and the proposed method (Degrees of freedom = 10,077,696).}
  \label{verifi:strong_speed-up_large}
\end{figure*}

Figures~\ref{verifi:strong_speed-up_small} and \ref{verifi:strong_speed-up_large} present the speed-up factors for the total FOM runtime, the total ROM runtime, and the two ROM components.
The total runtime of the proposed method exhibits favorable strong-scaling behavior in all analysis cases.
The reduced linear solver also scales well overall, although its scalability deteriorates at larger process counts.
For the case with 1,030,301 degrees of freedom shown in Fig.~\ref{verifi:strong_speed-up_small}, favorable solver scaling is maintained for $N_{\mathrm{POD}}=2{,}048$.
By comparison, the deterioration in solver speed-up becomes apparent at 64 processes for $N_{\mathrm{POD}}=512$ and at 128 processes for $N_{\mathrm{POD}}=1{,}024$.

This behavior can be attributed to the size of the reduced linear system.
A smaller value of $N_{\mathrm{POD}}$ results in a smaller reduced dimension and hence a lower computational workload per process.
As the process count increases, this workload becomes insufficient to amortize the communication and synchronization overheads, which consequently dominate the linear-solver runtime.
The larger FOM case shown in Fig.~\ref{verifi:strong_speed-up_large} exhibits a similar solver-scaling trend.
As indicated in Table~\ref{tab:values_strong1}, for a fixed number of POD computation subdomains, the reduced dimension is nearly independent of the number of full-order degrees of freedom.
The reduced linear systems therefore have comparable sizes in the two FOM cases, resulting in similar strong-scaling behavior of the linear solver.

\subsubsection{Verification of Static Load Balancing}\label{kensyo:load_balancing}

As discussed in Section~\ref{syou:staticlb}, the proposed method formulates static load balancing for domain-decomposed Galerkin ROMs as a problem of assigning weights to metanodes. To demonstrate the usefulness of this framework, we present an example of metanode weight assignment for static load balancing and quantitatively evaluate the parallel efficiency. Because the number of basis vectors stored in each POD computation subdomain differs from subdomain to subdomain, the computational load can become imbalanced among processes in parallel computations, potentially reducing the overall computational efficiency. To address this issue, this study introduces and verifies load balancing based on weighted graphs. Hyper-reduction is also introduced to further improve the computational efficiency. In this study, ECM (empirical cubature method) \cite{HERNANDEZ2017} is used as the hyper-reduction method, although the proposed method can be extended to arbitrary hyper-reduction methods \citep{CHATURANTABUT20102737,RYCKELYNCK2005346,WILLCOX2006208,KANEKO2021104385,Farhat2015,GRIMBERG20211846,CARLBERG2013623,LAUZON2024B474,HERNANDEZ2017}. The detailed formulation of the adopted hyper-reduction method is provided in Appendix B.

Let the set of all elements be $\mathcal{E}=\{e_1,\ldots,e_{n_{\mathrm{elem}}}\}$. Here, $n_{\mathrm{elem}}$ denotes the total number of elements. In a hyper-reduction method combined with domain decomposition, for each subdomain \(i\), a subset \(\tilde{\mathcal{E}}^{(i)}\subset\mathcal{E}\) consisting of $n^{(i)}_{\mathrm{HROM}}$ elements (\(\sum_{i=1}^{N_{\mathrm{POD}}} n^{(i)}_{\mathrm{HROM}}\ll n_{\mathrm{elem}}\)) and positive weights \(w_e >0\ (e\in\tilde{\mathcal{E}}^{(i)})\) are obtained, and the following approximations are used:
\begin{align}
  \bm{\Phi}^{\mathrm{T}} \bm{r}
  &
   \;\approx\;
   \sum_{i=1}^{N_{\mathrm{POD}}} \sum_{e\in\tilde{\mathcal{E}}^{(i)}} w_e\, \bm{\Phi}_e^{\mathrm{T}} \bm{r}_e, 
   \label{eq:ECM_r}\\
  \bm{\Phi}^{\mathrm{T}} \bm{J}\, \bm{\Phi}
  &
   \;\approx\;
   \sum_{i=1}^{N_{\mathrm{POD}}} \sum_{e\in\tilde{\mathcal{E}}^{(i)}} w_e\, \bm{\Phi}_e^{\mathrm{T}} \bm{J}_e \bm{\Phi}_e. 
   \label{eq:ECM_J}
\end{align}
Here, \(\bm{r}_e\in\mathbb{R}^{n_e}\) is the residual vector associated with finite element \(e\), \(\bm{J}_e\in\mathbb{R}^{n_e\times n_e}\) is the Jacobian associated with finite element \(e\), \(\bm{\Phi}_e\in\mathbb{R}^{n_e\times \sum_{i=1}^{N_{\mathrm{POD}}}n^{(i)}_{\mathrm{POD}}}\) is the POD basis assigned to finite element \(e\), and \(n_e\) is the number of degrees of freedom of finite element \(e\). The element set \(\tilde{\mathcal{E}}^{(i)}\) selected by the hyper-reduction method and the weights \(w_e\) are determined by solving an NNLS (non-negative least-squares) problem. In this study, the NNLS problem was formulated using matrices and vectors constructed from the residual, and the elements selected by this problem were identified. The details of the formulation are given in Appendix~B.2.

Thus, the computational cost of assembling the reduced system of linear equations in the ROM analysis depends on the number of elements selected in each POD computation subdomain, \(n^{(i)}_{\mathrm{HROM}}\). Therefore, in this study, the node weight \(n_{\mathrm{weight}}^{(i)}\) was set to either \(1\) or \(n^{(i)}_{\mathrm{HROM}}\), and the effect of load balancing based on the node weights was verified.

\begin{figure*}[t]
  \centering
  \begin{subfigure}[t]{0.4\linewidth}
    \hspace{25pt}
    \includegraphics[scale=0.5]{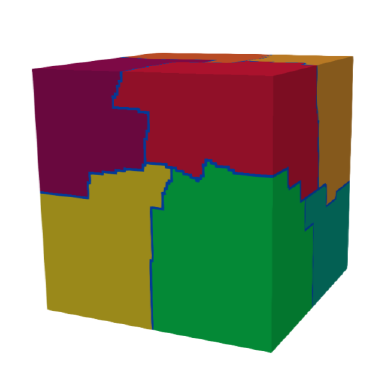}
    \subcaption{\small Before static load balancing ($n^{(i)}_{\mathrm{weight}}=1$).}
  \end{subfigure}
  \hspace{0.05\linewidth}
  \begin{subfigure}[t]{0.4\linewidth}
    \hspace{25pt}
    \includegraphics[scale=0.5]{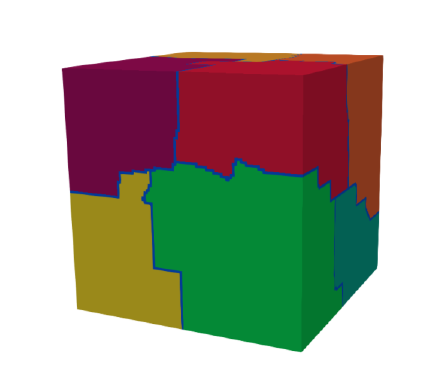}
    \subcaption{\small After static load balancing ($n^{(i)}_{\mathrm{weight}}=n^{(i)}_{\mathrm{HROM}}$).}
  \end{subfigure}
  \vspace{5.5pt} 
    \caption{\baselineskip=13pt Example of domain decomposition for POD computation subdomains before and after static load balancing
  ($N_{\mathrm{POD}}=512$, $N_{\mathrm{par}}=8$).}
    \vspace{-5pt}
    \label{fig:DDM_mesh_lb}
\end{figure*}

Unless otherwise stated, the basic computational settings were identical to those described in Section~\ref{kensyo:calceff}. 
First, finite element analyses were performed over the initial 100 time steps, and the snapshot matrix was constructed from the solution vectors obtained during this process. The snapshot matrix was then subjected to singular value decomposition to generate the POD basis used in the ROM analysis. The number of basis vectors used in each POD computation subdomain was determined according to Eq.~(\ref{eq:base_sel}). Next, the initial 100 time steps were analyzed using the ROM, and the element set $\tilde{\mathcal{E}}^{(i)}$ and weights $w_e\ (e\in\tilde{\mathcal{E}^{(i)}})$ were determined based on the hyper-reduction method. The number of POD computation subdomains was set to 512, and the number of parallel processes was varied as 8, 16, 32, 64, and 128. The number of compute nodes used was set to 1 for 8 and 16 parallel processes and to 2, 4, and 8 for 32, 64, and 128 parallel processes, respectively. The threshold for basis selection was set to $\varepsilon_{\mathrm{POD}}=10^{-6}$, and the threshold for element selection in the hyper-reduction method was set to $\varepsilon_{\mathrm{HROM}}=10^{-6}$. Figure~\ref{fig:DDM_mesh_lb} shows the domain decomposition before and after load balancing. For each POD computation subdomain $i$, the node weight $n^{(i)}_{\mathrm{weight}}$ was set to either 1 or $n^{(i)}_{\mathrm{HROM}}$.

\begin{figure}[t]
    \begin{center}
    \includegraphics[scale=0.44]{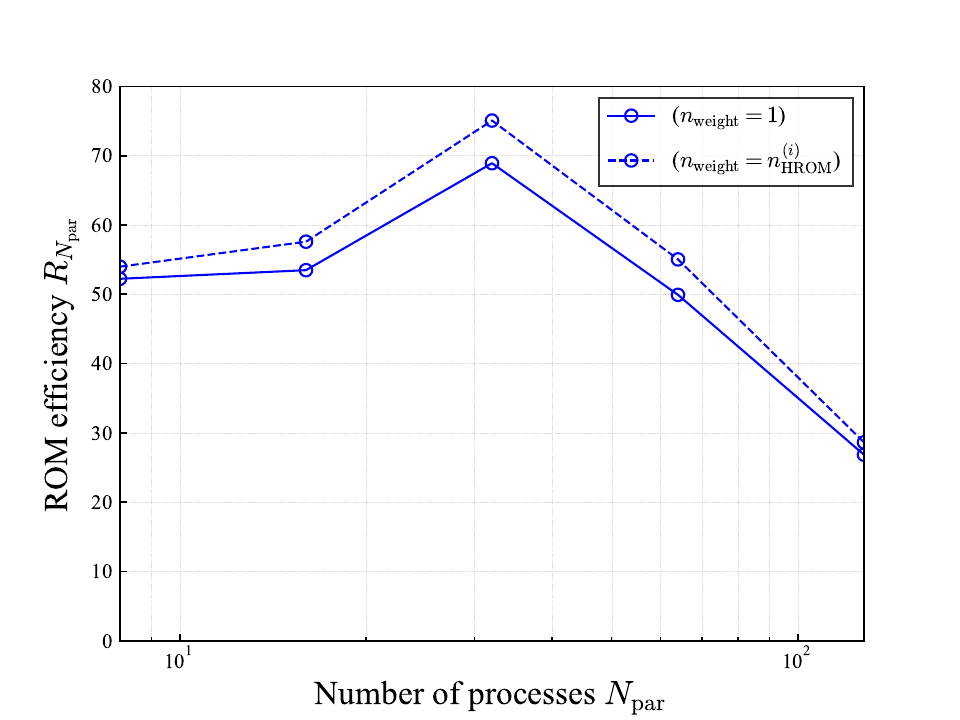}
    \end{center}
    \vspace{-10pt}
    \caption{\baselineskip=13pt ROM efficiency based on total runtime for different numbers of parallel processes.}
    \vspace{-5pt}
    \label{verifi:strong_ROM_efficiency_LB}
\end{figure}

\begin{figure}[t]
    \begin{center}
    \includegraphics[scale=0.44]{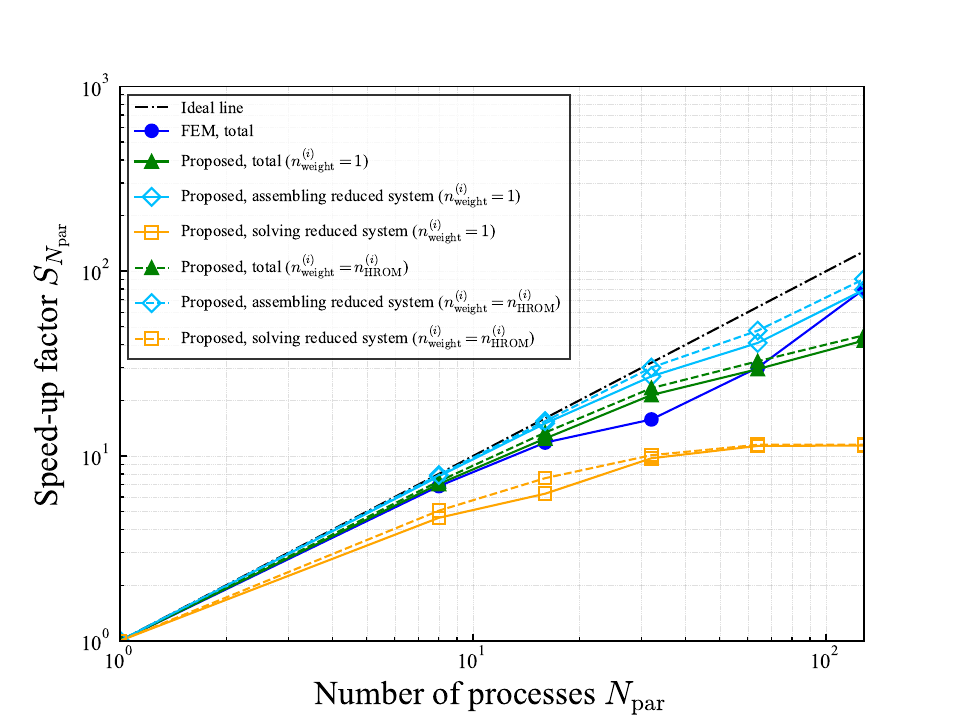}
    \end{center}
    \vspace{-10pt}
    \caption{\baselineskip=13pt Speed-up factors of the conventional FEM and the proposed method.}
    \vspace{-5pt}
    \label{verifi:strong_speed-up_small_LB}
\end{figure}

Figure~\ref{verifi:strong_ROM_efficiency_LB} shows the ROM efficiency for each number of parallel processes when the node weight $n_{\mathrm{weight}}^{(i)}$ in POD computation subdomain $i$ was set to either 1 or $n^{(i)}_{\mathrm{HROM}}$. The results show that setting the weight to $n^{(i)}_{\mathrm{HROM}}$ improved the computational efficiency by up to approximately 15\%. Figure~\ref{verifi:strong_speed-up_small_LB} shows the speed-up factors. The results show that the speed-up factor was consistently improved under all computational conditions. These results suggest that using $n^{(i)}_{\mathrm{HROM}}$ as the node weight is effective for improving the parallel efficiency and scaling behavior of the ROM. However, no marked improvement was observed for 8 and 128 parallel processes. On the other hand, the proposed method has the advantage that arbitrary node weights can be assigned; for example, further efficiency improvements are expected by optimizing the weights while accounting for imbalance in the number of basis vectors among subdomains. In ROM analyses involving parametric studies, the computational load for the next run can be estimated sequentially from previous execution results. Based on this estimate, the overall computation time can be reduced by allocating and optimizing the required resources and the number of hyper-reduction samples in advance. Dynamic load balancing and a generally applicable method for assigning weights remain topics for future work.

\subsection{Incompressible Navier--Stokes Equations}\label{sec:NSeq}
This section presents numerical simulations of unsteady flow around a three-dimensional cylinder as an application involving more complex physical phenomena and mesh structures than the numerical example in the previous section. This benchmark problem is widely used to evaluate the performance of various FOM and ROM fluid solvers, and has been adopted as a verification problem for ROMs in studies such as \citep{REYES2020112844,DAVE2025118194}.

In contrast to the linear diffusion problem considered in the previous section, the incompressible Navier--Stokes equations considered here are nonlinear. Accordingly, the Jacobian is updated at every Newton--Raphson iteration. In the runtime decomposition defined in the previous section, component~(a) therefore includes the assembly of both the reduced Jacobian and the reduced residual vector. For each time step, the reported runtime is accumulated over all Newton--Raphson iterations.


The momentum and continuity equations, which are the governing equations for an incompressible viscous fluid, are given in Eqs.~(\ref{eq:momentum1}) and (\ref{eq:continuity1}), respectively. The Dirichlet and Neumann boundary conditions are given in Eqs.~(\ref{eq:dirichlet1}) and (\ref{eq:neumann1}), respectively.

\begin{align}
\rho \left( \frac{\partial \bm{u}}{\partial t} 
  + \left(\bm{u} \cdot \nabla \right) \bm{u} \right)
   - \nabla \cdot \boldsymbol{\sigma} &= \bm{0}
  && \text{in } \Omega, \label{eq:momentum1} \\[6pt]
\nabla \cdot \bm{u} &= 0 
  && \text{in } \Omega, \label{eq:continuity1} \\[6pt]
\bm{u} &= \bm{u}_D 
  && \text{on } \Gamma_g, \label{eq:dirichlet1} \\[6pt]
\boldsymbol{\sigma}\bm{n} &= \bm{h} 
  && \text{on } \Gamma_h. \label{eq:neumann1}
\end{align}
Here, $\rho$ is the density, $\bm{u}$ is the velocity, and $p$ is the pressure. The Cauchy stress tensor $\boldsymbol{\sigma}$ is defined as
\begin{equation}
\boldsymbol{\sigma} 
= -p \bm{I} 
  + \mu \left( \nabla \bm{u} + \nabla \bm{u}^{\mathrm{T}} \right),
\end{equation}
where $\mu$ is the viscosity coefficient and $\bm{I}$ is the second-order identity tensor.

In this verification, common spatial and temporal discretization settings were used for the numerical method employed to generate the basis in the offline phase and for the numerical method to be reduced in the online phase. The finite element method was used for the spatial discretization, and the implicit Euler method was used for the temporal discretization. For the spatial discretization, a monolithic approach was adopted in which the velocity and pressure fields are treated as simultaneous unknowns. Galerkin least-squares stabilization \cite{TEZDUYAR19911,TEZDUYAR1992221} was introduced to suppress instabilities associated with advection and incompressibility. The detailed formulation of the discretization, the parameter settings, and the formulation of the nonlinear discrete equations based on Newton's method are described in Appendix A.2.

\begin{figure*}[t]
  \centering
  \begin{subfigure}[t]{0.38\linewidth}
    \hspace{-5pt}
    \includegraphics[scale=0.5]{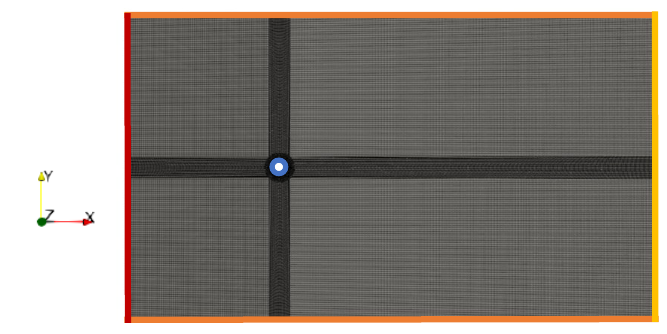}
    \subcaption{Computational mesh in the cross-section perpendicular to the $z$-axis.}
  \end{subfigure}
  \hspace{0.01\linewidth}
  \begin{subfigure}[t]{0.4\linewidth}
    \hspace{-10pt}
    \raisebox{15pt}{\includegraphics[scale=0.5]{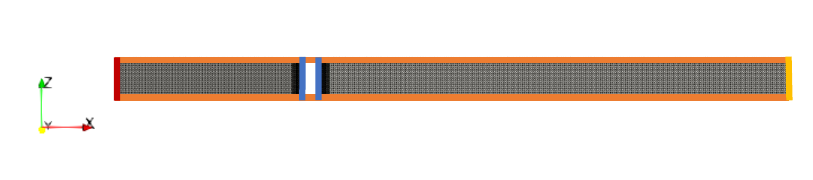}}
    \subcaption{Computational mesh in the cross-section perpendicular to the $y$-axis.}
  \end{subfigure}
  \vspace{5.5pt} 
  \caption{Computational mesh for flow past a three-dimensional cylinder.}
\label{flow_cylinder:mesh}
\end{figure*}

This verification was performed for a Reynolds number of \(\mathrm{Re}=200\). The representative length of the domain was set to the cylinder diameter \(D=1\,\mathrm{m}\), the density was set to \(\rho=200\,\mathrm{kg\,m^{-3}}\), and the viscosity coefficient was set to \(\mu=1.0\,\mathrm{Pa\cdot s}\). The computational mesh consists of hexahedral elements, with 2{,}488{,}800 elements, 2,629,200 nodes, and 10,516,800 total degrees of freedom. First-order Lagrange interpolation functions were used as the finite element basis functions. The mesh structure used in this study is shown in Fig.~\ref{flow_cylinder:mesh}. The time-step size $\Delta t$ was set to $0.005\,$s and was kept constant for both the FOM and ROM. The boundary conditions were assigned according to the colors in Fig.~\ref{flow_cylinder:mesh}(a) and (b): a uniform inflow velocity of \(U=1\,\mathrm{m/s}\) in the $x$ direction was imposed on the inflow boundary (red), a no-slip condition was imposed on the cylinder surface (cyan), slip conditions were imposed on the upper, lower, and spanwise boundaries (orange), and a traction-free condition was imposed on the outflow boundary (yellow). The origin of the coordinate system was placed at the cylinder center; the inflow and outflow boundaries were located at \(x=-10D\) and \(x=25D\), respectively, the upper and lower boundaries were located at $y=\pm 10D$, and the spanwise extent was $z\in[-D,D]$

\begin{figure*}[t]
  \centering
  \begin{subfigure}[t]{0.4\linewidth}
    \includegraphics[scale=0.4]{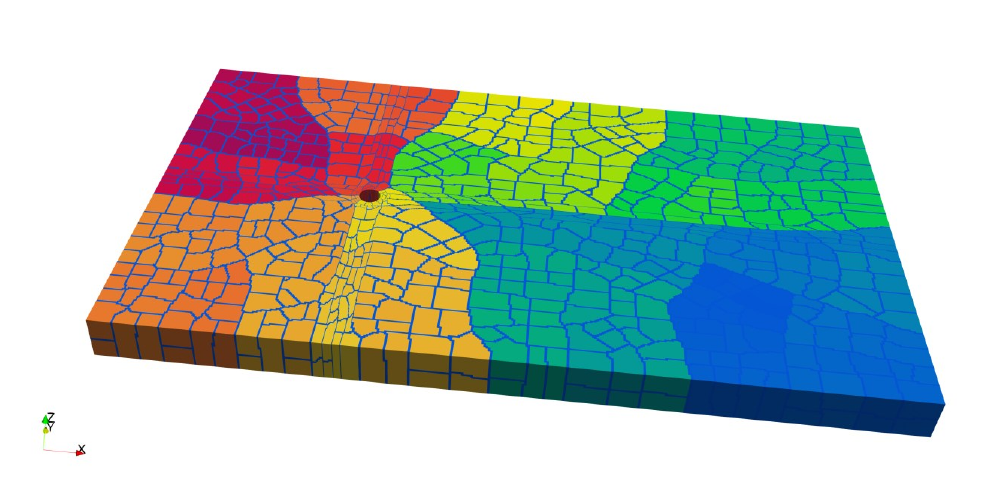}
    \subcaption{POD computation subdomains ($N_{\mathrm{POD}}=1024$).}
  \end{subfigure}
  \hspace{0.05\linewidth}
    \begin{subfigure}[t]{0.4\linewidth}
    \includegraphics[scale=0.4]{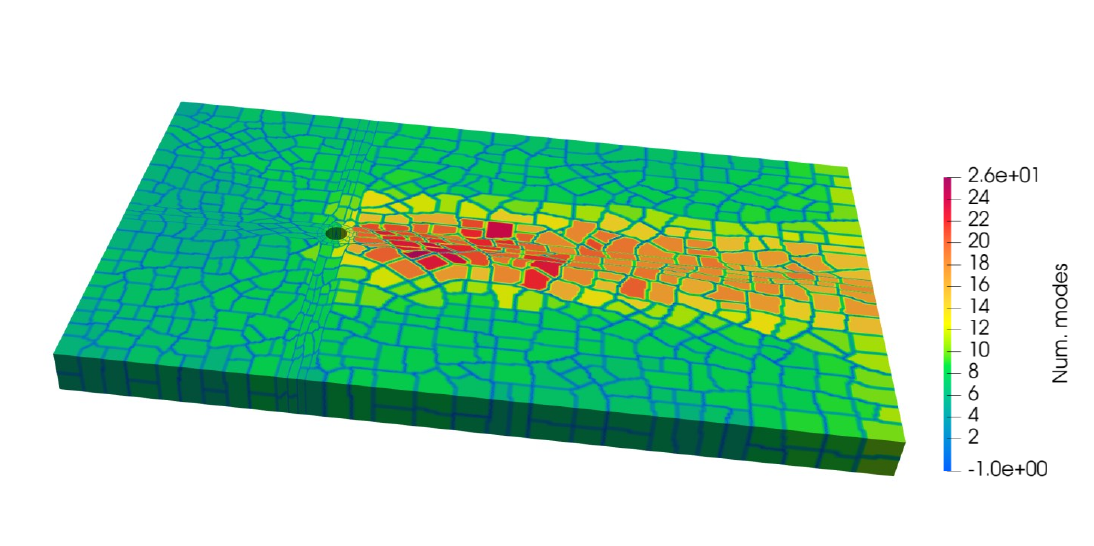}
    \subcaption{Number of basis vectors assigned to each subdomain (overlapping elements are shown in blue).}
  \end{subfigure}
  \begin{subfigure}[t]{0.4\linewidth}
    \includegraphics[scale=0.4]{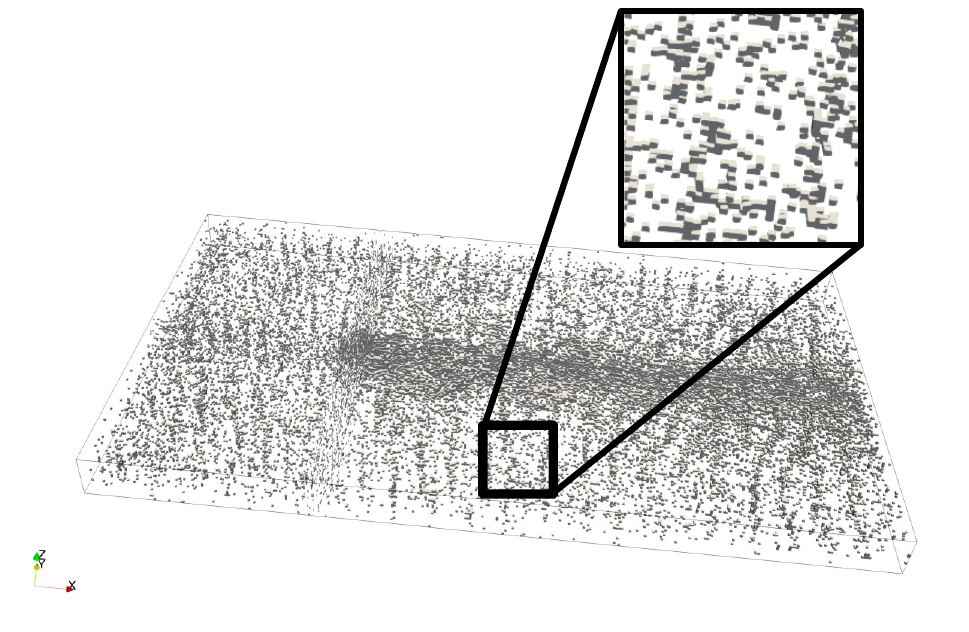}
    \subcaption{Finite elements selected by hyper-reduction (sampled elements are shown in gray).}
  \end{subfigure}
    \vspace{5.5pt} 
  \caption{Domain decomposition for local POD computations and hyper-reduction sampling.}
\label{Fig:ECMmesh_cylinder}
\end{figure*}

The purpose of this numerical example is to reproduce the periodic vortex shedding in the cylinder wake using the proposed ROM framework. After the flow had sufficiently developed and periodic vortex shedding appeared steadily, the flow field at \(t=100\,\mathrm{s}\) was used as the initial condition, and snapshots were collected for POD over a 1\,s interval, \(100 \le t \le 101\,\mathrm{s}\). Snapshots were collected at a sampling time interval of \(\Delta t_{\mathrm{snap}}=0.01\,\mathrm{s}\), and a total of 101 snapshots were used. The ROM analysis was performed over the same interval \([100,101]\,\mathrm{s}\), and model reduction was carried out using the obtained POD basis. Using the ROM solution over the interval \([100,101]\,\mathrm{s}\), the element set $\tilde{\mathcal{E}}$ and weights $w_e,(e\in\tilde{\mathcal{E}})$ were determined based on the hyper-reduction method. The number of POD computation subdomains was set to 1024, and the number of parallel processes was varied as 32, 64, 128, 256, and 512. The number of compute nodes used for 32, 64, 128, 256, and 512 parallel processes was set to 2, 2, 4, 8, and 16, respectively. The threshold for basis selection was set to $\varepsilon_{\mathrm{POD}}=10^{-6}$, and the threshold for element selection in the hyper-reduction method was set to $\varepsilon_{\mathrm{HROM}}=10^{-6}$. The BiCGSafe method \cite{FUJINO2005145} was used to solve the systems of linear equations, with a convergence tolerance of $1.0\times 10^{-8}$ and a maximum of 10,000 iterations. Diagonal scaling was used as the preconditioner. Figure~\ref{Fig:ECMmesh_cylinder} shows the visualization of the POD computation subdomains, the distribution of the number of basis vectors in each POD computation subdomain, and the elements selected by the hyper-reduction method.

\begin{figure}[t]
    \begin{center}
    \hspace{25pt}
    \includegraphics[scale=0.6]{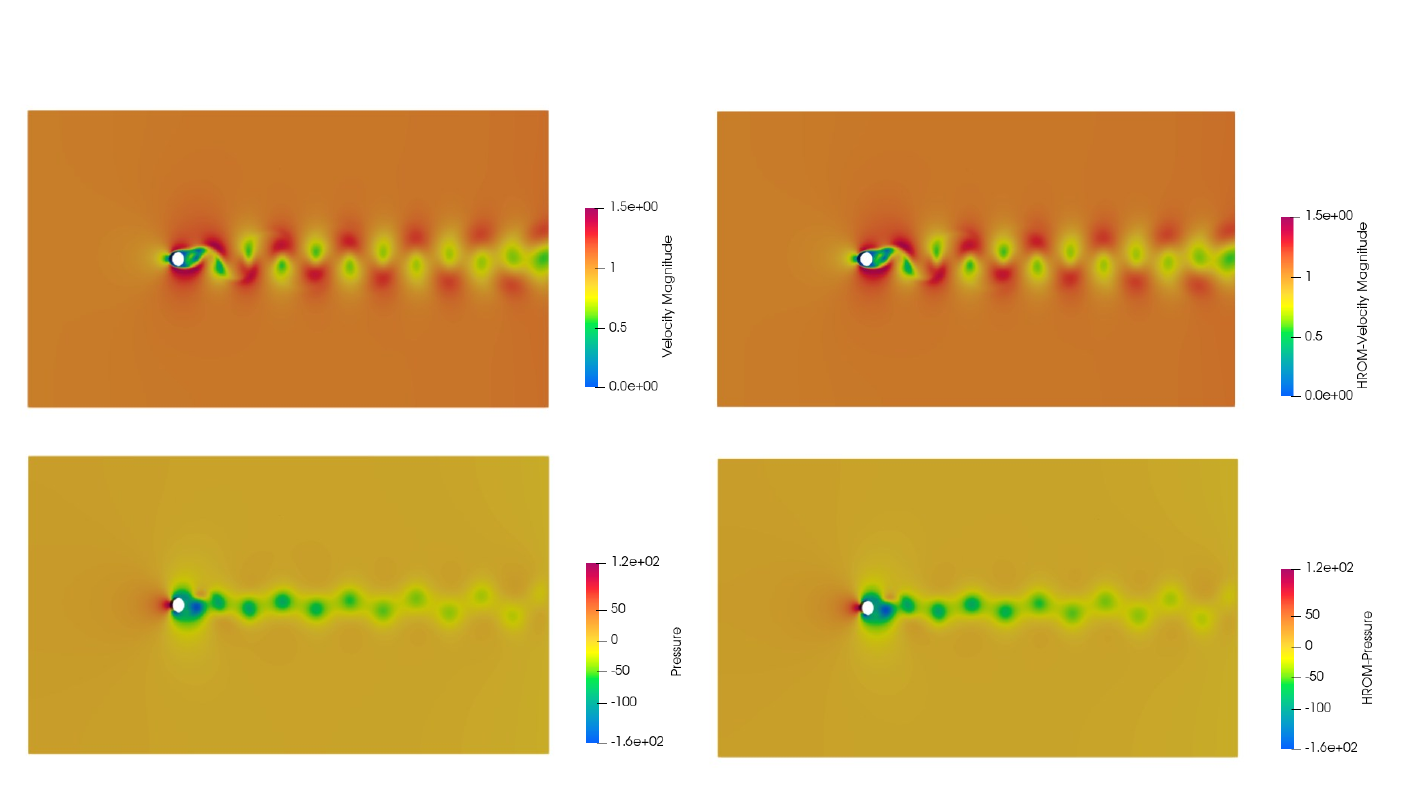}
    \end{center}
    \vspace{-10pt}
    \caption{\baselineskip=13pt Contours of velocity magnitude (top) and pressure (bottom) at t = 101 s for Re = 200 on a cross-section normal to the $z$-axis: FOM (left column) and ROM (right column).}
    \vspace{-5pt}
    \label{Fig:contours_cylinder}
\end{figure}

\begin{figure}[t]
    \begin{center}
    \includegraphics[scale=0.44]{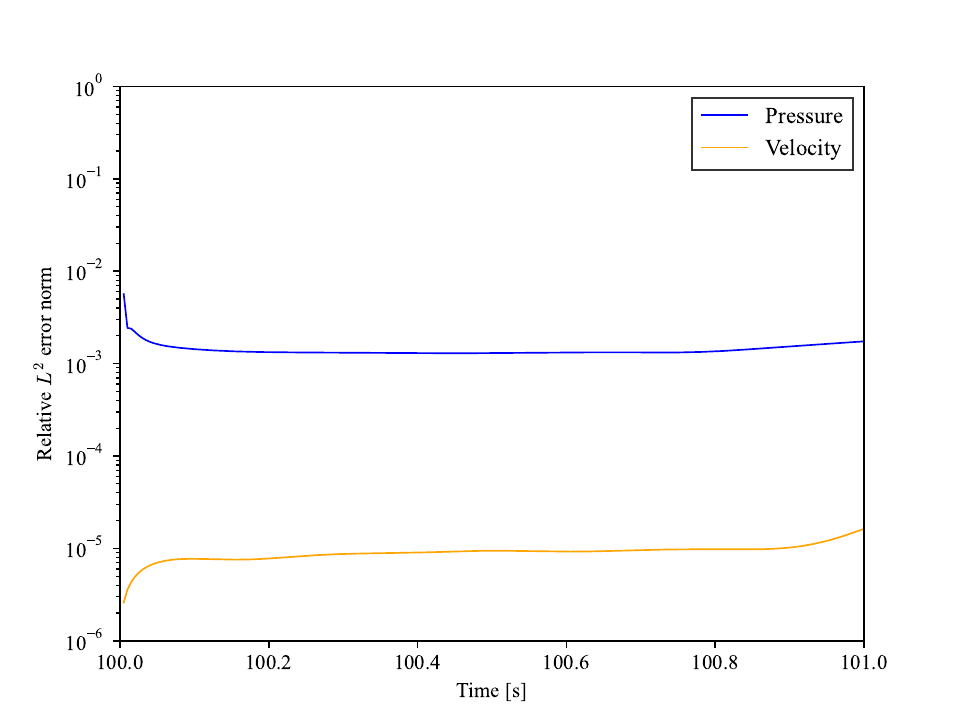}
    \end{center}
    \vspace{-10pt}
    \caption{\baselineskip=13pt Relative $L^2$ error norms of velocity and pressure between FOM and ROM using ECM.}
    \vspace{-5pt}
    \label{verifi:l2error_cylinder}
\end{figure}

Figure~\ref{Fig:contours_cylinder} shows the ROM visualization results for the velocity and pressure fields for the case with the basis selection threshold \(\varepsilon_{\mathrm{POD}}=10^{-6}\) and the element selection threshold \(\varepsilon_{\mathrm{HROM}}=10^{-6}\). The results show that the flow field obtained by the FOM is qualitatively reproduced by the ROM. Figure~\ref{verifi:l2error_cylinder} shows the computational accuracy of the velocity and pressure fields for the same case. The relative error of the velocity field remained approximately on the order of $10^{-5}$, whereas that of the pressure field remained approximately on the order of $10^{-3}$. Although the error in the pressure field tended to be larger than that in the velocity field, it was also kept at an error level of approximately $10^{-3}$, confirming that the present ROM achieves good approximation accuracy for the pressure field as well.

\begin{figure}[t]
    \begin{center}
    \includegraphics[scale=0.44]{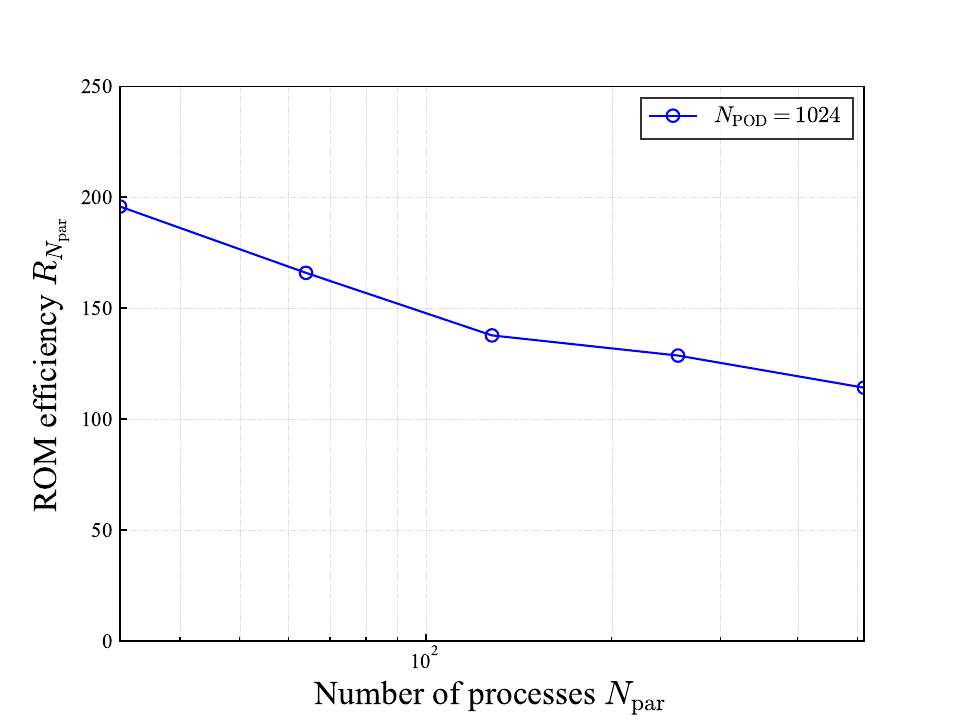}
    \end{center}
    \vspace{-10pt}
    \caption{\baselineskip=13pt ROM efficiency based on total runtime for different numbers of parallel processes.}
    \vspace{-5pt}
    \label{verifi:strong_ROM_efficiency_cylinder}
\end{figure}

\begin{figure}[t]
    \begin{center}
    \includegraphics[scale=0.44]{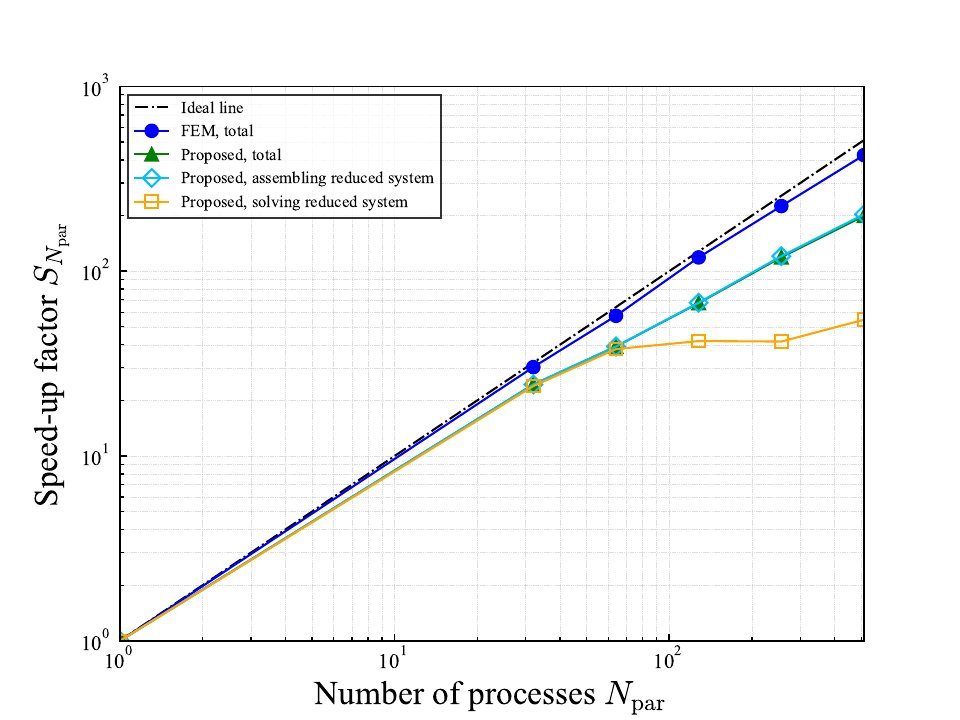}
    \end{center}
    \vspace{-10pt}
    \caption{\baselineskip=13pt Speed-up factors of the conventional FEM and the proposed method.}
    \vspace{-5pt}
    \label{verifi:strong_sprrdupfactor_cylinder}
\end{figure}

The ROM efficiency defined in Eq.~(\ref{eq:ROM_efficiency}) is shown in Fig.~\ref{verifi:strong_ROM_efficiency_cylinder}. In each analysis case, the computation time was evaluated over the first 10 time steps, and the average of three measurements was used for evaluation. The results show that a speedup of approximately 100--200 times relative to the FOM was obtained for all numbers of parallel processes, demonstrating that the proposed method also achieves high parallel computing performance in the analysis of flow around a cylinder. Figure~\ref{verifi:strong_sprrdupfactor_cylinder} shows the speed-up factors. The results show that the speed-up factor of the total computation time of the proposed method remains high. The speed-up factor of the linear solver for the reduced system of linear equations is generally good, although it decreases as the number of parallel processes increases. As in Sections~\ref{kensyo:calceff} and \ref{kensyo:load_balancing}, this is likely because, as the number of parallel processes increases, the computation time within each process in the linear solver becomes short, making communication and synchronization overheads dominant relative to the computation time.

%% file: include/conclusions.tex
\clearpage
\section{Conclusions}
This study proposed the metagraph-based domain-decomposed Galerkin reduced-order model (MBDD-G-ROM), which combines a graph-based generalization of DD-G-ROMs with a hierarchical domain-partitioning strategy for distributed-memory parallel computation. The proposed method represents DD-G-ROMs defined over arbitrary domain decompositions using calculation-point graphs and metagraphs, thereby providing a unified description of subdomain connectivity and the associated block structure of the reduced system. Based on this representation, the hierarchical domain-partitioning strategy decouples the POD computation subdomains from the parallel computation subdomains, enabling distributed-memory parallelization of L-POD, including both reduced-system assembly and the linear solver.

More specifically, each POD computation subdomain is represented as a metanode, and metaedges are defined according to the overlap of the supports of the local POD basis functions associated with different subdomains. The resulting metagraph encodes the block-sparsity pattern of the L-POD reduced system. This representation makes it possible to treat arbitrary-shaped POD computation subdomains as metanodes and to assign them independently to parallel computation subdomains through the second-stage graph partitioning. Moreover, this metagraph structure provides a natural way to incorporate static load balancing through user-defined metanode weights representing estimated computational costs.

The proposed method was then evaluated in terms of ROM accuracy and parallel performance. In the verification based on the unsteady diffusion equation, increasing the number of POD computation subdomains improved the ROM accuracy, demonstrating the expected behavior of the L-POD implementation within the proposed method. Strong-scaling tests showed that the proposed method achieved high parallel efficiency even for large-scale problems. The static load-balancing test further demonstrated that metanode weights based on estimated computational costs allow the second-stage graph partitioning to account for cost differences among subdomains. In the tested case, this weighting improved computational efficiency. Finally, the proposed method was applied to flow around a cylinder, demonstrating its applicability to complex flow fields and its ability to maintain high parallel computing efficiency.

Although this study employed L-POD to construct local linear reduced bases, the proposed metagraph-based generalization is not restricted to POD-based basis construction. It can, in principle, accommodate a broad class of local model-reduction techniques whose inter-subdomain interactions can be represented on the metagraph. Future work will therefore include extending the present POD-based linear-subspace implementation to incorporate nonlinear-manifold-based approaches \citep{RUTZMOSER2017196,BARNETT2022111348,GEELEN2023115717,KIM2024116978,LEE2020108973,DIAZ2024116943}, with the aim of improving ROM accuracy and efficiency for problems with strongly nonlinear solution structures. In such problems, the effective dimension of the local low-dimensional representations and the associated computational cost may vary with time and parameters. Dynamic load balancing within the proposed metagraph-based framework will also be investigated.

\newpage

%% file: include/appendixA.tex
\newpage
\appendix
\renewcommand{\thesection}{\Alph{section}}
\renewcommand{\thesubsection}{\Alph{section}.\arabic{subsection}}
\newpage

\section{Finite Element Method}

In this appendix, the analysis domain is assumed to be a three-dimensional bounded domain $\Omega\subset\mathbb{R}^3$, and its boundary is denoted by $\Gamma$. The boundary $\Gamma$ is divided into the Dirichlet boundary $\Gamma_g$ and the Neumann boundary $\Gamma_h$, satisfying $\Gamma=\Gamma_g\cup\Gamma_h$ and $\Gamma_g\cap\Gamma_h=\emptyset$.

\subsection{Unsteady Diffusion Equation}
The unsteady diffusion equation used as the governing equation is given by
\begin{equation}\label{eq:poisson_equ_appendix}
    \frac{\partial u}{\partial t}-\,\nabla\cdot ( k\,\nabla u)=Q,
\end{equation}
where $k$ is the diffusion coefficient, $u$ is the unknown representing the diffusing physical quantity, and $Q$ is the source term representing the physical quantity generated in the analysis domain $\Omega$. The boundary conditions are given by
\begin{align}
    u&=u_D\quad \text{on } \Gamma_g, \label{eq:dirichlet_bc_A}\\
    \nabla {u}\cdot \bm{n}&= {q}_N \quad \text{on } \Gamma_h, \label{eq:neumann_bc_A}
\end{align}
where $u_D$ is the value of the solution $u$ on the Dirichlet boundary $\mathrm{\Gamma}_g$, $\bm{n}$ is the outward normal vector on the Neumann boundary $\mathrm{\Gamma}_h$, and $q_N$ is a function representing the physical quantity prescribed on the boundary $\mathrm{\Gamma}_h$.

The solution space and test function space are defined as
\begin{equation}
\mathcal{S}_h:=\{u\in H^1(\Omega)\mid u|_{\Gamma_g}=u_D\}, \qquad
\mathcal{V}_h:=\{w\in H^1(\Omega)\mid w|_{\Gamma_g}=0\}.
\end{equation}

The weak form of this problem is given as follows:

\begin{equation}\label{diff_weakfoam}
\exists\,u\in\mathcal{S}_h\ \text{s.t.}\ \forall w\in\mathcal{V}_h:\ a_\Omega(u,w)=\ell_\Omega(w),
\end{equation}
where the bilinear form $a_\Omega(u,w)$ and the linear functional $\ell_\Omega(w)$ are defined by
\begin{equation}
  a_\Omega(u,w) = \Bigl(\frac{\partial u}{\partial t},\,w\Bigr)_\Omega 
                  +\,(k\nabla u,\,\nabla w)_\Omega,
  \qquad
  \ell_\Omega(w) = (Q,\,w)_\Omega + \langle k\,q_N,\,w\rangle_{\Gamma_h}.
\end{equation}

Using the basis functions $N_i$ of the finite element space $\mathcal{S}_h$, the unknown $u$ is expressed as
\begin{equation}\label{eq:fem_solv_vec}
    u\approx \sum_{i=1}^{n_\mathrm{FOM}} N_i\, u_i, 
\end{equation}
where $n_\mathrm{FOM}$ is the number of finite element nodes, $N_i$ is the basis function defined at finite element node $i$, and $u_i$ is the coefficient corresponding to $N_i$. The time derivative term is approximated by the backward difference formula
\begin{equation}\label{eq:euler_implicit_time}
    \frac{\partial u_i}{\partial t}
    \approx
    \frac{u_i^{t+\Delta t}-u_i^{t}}{\Delta t},
\end{equation}
where $u_i^{t}$ is the coefficient corresponding to the basis function $N_i$ at time $t$, and $\Delta t$ is the time-step size. In this study, the diffusion and source terms are evaluated at time $t+\Delta t$, yielding a temporal discretization based on the implicit Euler method. By applying Eqs.~\eqref{eq:fem_solv_vec} and \eqref{eq:euler_implicit_time}, together with the Galerkin method, to the weak form \eqref{diff_weakfoam}, the following discretized equation is obtained:
\begin{equation}\label{eq:weight_res}
    \bm{w}^{\mathrm{T}}\bm{K}\bm{u}=\bm{w}^{\mathrm{T}}\bm{f},
\end{equation}
where the matrix and vector components are given by
\begin{align}
    K_{ij}&= \int_{\Omega}{ N_i N_j \,\mathrm d\Omega} + \Delta t \int_{\Omega}{k\,\nabla N_i \cdot \nabla N_j \,\mathrm d\Omega},\\
    u_{i}&= u^{t+\Delta t}_i,\\
    f_{i}&=\Delta t\int_{\Omega}{N_i\, Q \,\mathrm d\Omega}
    +\Delta t\int_{\Gamma_h}{k\,N_i\, q_N \,\mathrm d\Gamma}
    +  \sum_{j=1}^{n_\mathrm{FOM}} \left(\int_{\Omega}{ N_i N_j \,\mathrm d\Omega}\right)\, u^{t}_j.
\end{align}
\hspace{-1mm}Here, $\bm{w}\in \mathbb{R}^{n_\mathrm{FOM}}$ is an arbitrary vector derived from the test function. The vector components $w_i$ of the weight vector $\bm{w}$ satisfy $w_i=0\ (i \in \eta_g)$, where $\eta_g$ denotes the set of indices of finite element nodes on which Dirichlet boundary conditions are imposed. The matrix $\bm{K}\in \mathbb{R}^{n_\mathrm{FOM} \times n_\mathrm{FOM}}$ is the coefficient matrix in the finite element method, $\bm{u}\in \mathbb{R}^{n_\mathrm{FOM}}$ is the solution vector, and $\bm{f}\in \mathbb{R}^{n_\mathrm{FOM}}$ is the right-hand-side vector. The solution vector $\bm{u}$ is expressed as $\bm{u}=\bm{\bar{u}}+\bm{g}$, where $\bm{\bar{u}}$ is the unknown vector and $\bm{g}$ is a known vector satisfying the Dirichlet boundary conditions. The vector components $\bar{u}_i$ of $\bm{\bar{u}}$ satisfy $\bar{u}_i=0\ (i \in \eta_g)$, and the vector components $g_i$ of $\bm{g}$ satisfy $g_i=0\ (i \in \eta-\eta_g)$. Here, $\eta$ is the index set of the finite element nodes. The function $\bm{g}$ prescribed to satisfy the Dirichlet boundary conditions is also referred to as a lifting function \cite{BALLARIN2015,OXBERRY2017}.

\subsection{Incompressible Navier--Stokes Equations}

The momentum and continuity equations governing an incompressible viscous fluid are given in Eqs.~\eqref{eq:momentum} and \eqref{eq:continuity}, respectively, and the Dirichlet and Neumann boundary conditions are given in Eqs.~\eqref{eq:dirichlet} and \eqref{eq:neumann}, respectively.
\begin{align}
\rho \left( \frac{\partial \bm{u}}{\partial t} 
  + \left(\bm{u} \cdot \nabla \right) \bm{u} - \bm f\right)
   - \nabla \cdot \bm{\sigma} &= \bm{0}
  && \text{in } \Omega, \label{eq:momentum} \\
\nabla \cdot \bm{u} &= 0 
  && \text{in } \Omega, \label{eq:continuity} \\
\bm{u} &= \bm{u}_D 
  && \text{on } \Gamma_g, \label{eq:dirichlet} \\
\bm{\sigma}\,\bm{n} &= \bm{h} 
  && \text{on } \Gamma_h, \label{eq:neumann}
\end{align}
where $\rho$ is the density, $\bm{u}$ is the velocity, $p$ is the pressure, $\bm f$ is the body force per unit mass, and $\bm n$ is the outward unit normal vector. The Cauchy stress tensor $\bm{\sigma}$ is given by
\begin{equation}
\bm{\sigma} 
= -p \bm{I} 
  + \mu \left( \nabla \bm{u} + (\nabla \bm{u})^{\mathrm{T}} \right),
\end{equation}
where $\mu$ is the viscosity coefficient and $\bm I$ is the second-order identity tensor.

For the spatial discretization, this study adopts a monolithic formulation in which the velocity and pressure fields are treated simultaneously as unknowns. To suppress instabilities associated with advection and incompressibility, Galerkin least-squares (GLS) stabilization \cite{TEZDUYAR19911,TEZDUYAR1992221} is introduced. Specifically, the streamline-upwind Petrov--Galerkin (SUPG), pressure-stabilizing Petrov--Galerkin (PSPG), and least-squares on incompressibility constraint (LSIC) terms are added.

First, the solution and test function spaces are defined as
\begin{equation}
  \mathcal{S}_h:=\bigl\{(\bm u_h,p_h)\in \mathcal{V}_h\times \mathcal{Q}_h\ \big|\ \bm u_h|_{\Gamma_g}=\bm u_D\bigr\},\quad
  \mathcal{W}_h:=\bigl\{(\bm w_h,q_h)\in \mathcal{V}_h\times \mathcal{Q}_h\ \big|\ \bm w_h|_{\Gamma_g}=\bm 0\bigr\}.
\end{equation}
Here, $\bm u_h$ is the finite element approximation of the velocity, $p_h$ is the finite element approximation of the pressure, $\bm w_h$ is the test function corresponding to the velocity, and $q_h$ is the test function corresponding to the pressure. Moreover, $\mathcal{V}_h\subset [H^1(\Omega)]^3$ and $\mathcal{Q}_h\subset L^2(\Omega)$.

Next, the residuals of the momentum and continuity equations are defined as
\begin{align}
\bm r_M(\bm u_h,p_h)&=
\rho \left(\frac{\partial \bm u_h}{\partial t} + (\bm u_h\!\cdot\!\nabla)\bm u_h -\bm f\right)
- \nabla\!\cdot\!\bm{\sigma}(\bm u_h,p_h), \\
r_C(\bm u_h)&=\nabla\!\cdot\!\bm u_h.
\end{align}

The weak form with GLS stabilization is given by

\begin{equation}
\exists\,(\bm u_h,p_h)\in\mathcal{S}_h\ \text{s.t.}\ \forall\,(\bm w_h,q_h)\in\mathcal{W}_h:\ 
B(\bm w_h,q_h;\bm u_h,p_h)=F(\bm w_h,q_h),
\end{equation}
where $B$ and $F$ are defined as
\begin{align}
B(\bm w_h,q_h;\bm u_h,p_h)
&= \int_{\Omega} \bm w_h \cdot
\rho\Big(\frac{\partial \bm u_h}{\partial t} 
+ (\bm u_h\!\cdot\!\nabla)\bm u_h\Big)\,\mathrm d\Omega \notag\\
&\quad + \int_{\Omega} \tfrac12\!\left(\nabla \bm w_h+(\nabla \bm w_h)^{\mathrm T}\right) : \bm{\sigma}(\bm u_h,p_h)\,\mathrm d\Omega
+ \int_{\Omega} q_h\, \nabla\!\cdot\!\bm u_h\, \mathrm d\Omega \notag\\
&\quad + \sum_{e=1}^{n_{\mathrm{elem}}} \int_{\Omega_e} \tau_M^e\,
\Big( (\bm u_h\!\cdot\!\nabla)\bm w_h + \nabla \frac{q_h}{\rho} \Big)\cdot
\bm r_M(\bm u_h,p_h)\, \mathrm d\Omega \notag\\
&\quad + \sum_{e=1}^{n_{\mathrm{elem}}}\int_{\Omega_e}\rho\,\tau_C^e\,(\nabla\!\cdot\!\bm w_h)\, r_C(\bm u_h)\,\mathrm d\Omega,\\
F(\bm w_h,q_h)&= \int_{\Omega} \bm w_h \cdot
\rho\, \bm f\, \mathrm d \Omega +\int_{\Gamma_h}  \bm w_h\!\cdot\!\bm h \,\mathrm d\Gamma,
\end{align}
where $n_{\mathrm{elem}}$ is the total number of finite elements. The element-wise stabilization parameters $\tau_M^e$ and $\tau_C^e$ are defined based on Taylor et al.~\cite{TAYLOR1998155} and Whiting et al.~\cite{WHITING200193} as follows:
\begin{align}
\tau_M^e &=
\left(
    \frac{4}{\Delta t^2}
    + \bm{u}_h\cdot\bm{G}\bm{u}_h
    + C_I \,\nu^2 \, \bm{G} : \bm{G}
\right)^{-1/2},\\
\tau_C^e &=\left( \operatorname{tr}(\bm{G}) \, \tau_M^e \right)^{-1},
\end{align}
where $\Delta t$ is the time-step size, $\nu=\mu/\rho$ is the kinematic viscosity, and $C_I$ is an element-dependent constant. In this study, $C_I =36$ is used, which is the value employed for first-order hexahedral elements \cite{WHITING1999}. The tensor $\bm{G}$ is the element metric tensor, defined as
\begin{equation}
\bm{G} = 
\left( \frac{\partial \bm{\xi}}{\partial \bm{x}} \right)^{\mathrm{T}}
\left( \frac{\partial \bm{\xi}}{\partial \bm{x}} \right),
\end{equation}
where $\bm x$ is the coordinate vector in the physical coordinate system and $\bm \xi$ is the coordinate vector in the natural coordinate system.

In this study, the weak form with the above GLS (SUPG/PSPG/LSIC) stabilization is spatially discretized using the finite element method, and the implicit Euler method is used for the temporal discretization. The resulting nonlinear system of equations is solved by Newton--Raphson iteration. The same function space is used for both velocity and pressure.

First, the velocity, pressure, and test functions are approximated using the basis functions $N_A$ as follows:
\begin{align}\label{appendixA_NSeq_eq:fem_solv_vec}
  \bm u_h &= \sum_{A=1}^{n_\mathrm{FOM}} N_A\,\bm u_A, &
  p_h &= \sum_{A=1}^{n_\mathrm{FOM}} N_A\,p_A, \\
  \bm w_h &= \sum_{A=1}^{n_\mathrm{FOM}} N_A\,\bm w_A, &
  q_h &= \sum_{A=1}^{n_\mathrm{FOM}} N_A\,q_A,
\end{align}
where $n_\mathrm{FOM}$ is the number of finite element nodes. The time derivative term is approximated by the backward difference formula
\begin{equation}\label{appendixA_NSeq_eq:euler_implicit_time}
    \frac{\partial \bm u_A}{\partial t}=\frac{\bm u_A^{t+\Delta t}-\bm u_A^{t}}{\Delta t},
\end{equation}
where $\bm u_A^{t}$ is the velocity vector corresponding to the basis function $N_A$ at time $t$. In this study, the convective, viscous, and pressure terms are evaluated at time $t+\Delta t$, yielding a temporal discretization based on the implicit Euler method.

Next, the unknown vector $\bm U_{t+\Delta t}$ at time $t+\Delta t$ is defined as
\begin{equation}
  \bm U_{t+\Delta t}
  =
  \begin{Bmatrix}
    \bm u^{t+\Delta t} \\
    \bm p^{t+\Delta t}
  \end{Bmatrix},\qquad
    \bm u^{t+\Delta t}=\begin{bmatrix}
    (\bm u_1^{t+\Delta t})^{\mathrm T}, &
    \cdots, &
    (\bm u_{n_\mathrm{FOM}}^{t+\Delta t})^{\mathrm T}
  \end{bmatrix}^{\mathrm T},\qquad
    \bm p^{t+\Delta t}=\begin{bmatrix}p_1^{t+\Delta t},&\cdots,&p_{n_\mathrm{FOM}}^{t+\Delta t}\end{bmatrix}^\mathrm{T}.
\end{equation}

After weighting by the test functions, the Newton--Raphson correction $\Delta\bm U^{\,k}_{t+\Delta t}$ satisfies
\begin{align}
\label{eq:block_system_imp_euler_residual}
  \bm W^{\mathrm T}\bm J\!\left(\bm U^{\,k}_{t+\Delta t}\right)\Delta\bm U^{\,k}_{t+\Delta t}
  &= -\,\bm W^{\mathrm T}\bm r\!\left(\bm U^{\,k}_{t+\Delta t}\right),\\
  \label{eq:update_NR_residual}
  \bm U^{\,k+1}_{t+\Delta t}
  &= \bm U^{\,k}_{t+\Delta t}
  + \Delta\bm U^{\,k}_{t+\Delta t},
\end{align}
where $k$ is the iteration index in the Newton--Raphson method.

The iteration is continued until the following condition is satisfied:
\begin{equation}
  \frac{
    \left\|\bm r\!\left(\bm U^{\,k}_{t+\Delta t}\right)\right\|
  }{
    \left\|\bm r\!\left(\bm U^{\,0}_{t+\Delta t}\right)\right\|
  }
  < \varepsilon_{\mathrm{NR}}, 
\end{equation}
where $\bm r\!\left(\bm U^{\,0}_{t+\Delta t}\right)$ is the initial residual vector, and $\varepsilon_{\mathrm{NR}}$ is the convergence tolerance for the Newton--Raphson method. In this study, $\varepsilon_{\mathrm{NR}}=10^{-6}$ is used.

The Jacobian $\bm J$, correction $\Delta\bm U^{\,k}_{t+\Delta t}$, residual vector $\bm r$, and arbitrary vector $\bm W$ derived from the test functions are written in block form as
\begin{equation}\label{eq:J_blocks}
\begin{aligned}
  \bm J=&
  \begin{bmatrix}
    \bm J_{MM} & \bm J_{MC} \\
    \bm J_{CM} & \bm J_{CC}
  \end{bmatrix},
  \qquad
  \Delta\bm U^{\,k}_{t+\Delta t} :=
  \begin{Bmatrix}
    \Delta\bm u^{\,k}_{t+\Delta t} \\
    \Delta \bm p^{\,k}_{t+\Delta t}
  \end{Bmatrix},
  \qquad
  \bm r=
  \begin{Bmatrix}
    \bm R_M \\
    \bm R_C
  \end{Bmatrix},
  \\
    \bm W
  =&
  \begin{Bmatrix}
    \bm w \\
    \bm q
  \end{Bmatrix},
  \qquad
\bm w=
\begin{bmatrix}
(\bm w_1)^{\mathrm T} ,& \cdots ,& (\bm w_{n_\mathrm{FOM}})^{\mathrm T}
\end{bmatrix}^{\mathrm T},
\qquad
\bm q=
\begin{bmatrix}
q_1, & \cdots, & q_{n_\mathrm{FOM}}
\end{bmatrix}^{\mathrm T}.
\end{aligned}
\end{equation}

Here, $\bm R_M$ and $\bm R_C$ are the residual vectors for the momentum and continuity equations, respectively. The components of the residual vectors $\bm R_M$ and $\bm R_C$ are defined using the basis and test functions as follows:
\begin{align}
  \bm R_M&=[(R_M)_{A,i}],\\
  (R_M)_{A,i}
  &= B(N_A\bm e_i,\,0;\,\bm u_h,p_h)\;-\;F(N_A \bm e_i,\,0),\\
  \bm R_C&=[(R_C)_A],\\
  (R_C)_A
  &= B(\bm 0,\,N_A;\,\bm u_h,p_h)\;-\;F(\bm 0,\,N_A).
\end{align}
Here, $\bm e_i\ (i=1,2,3)$ denotes an orthonormal basis vector. The above equations correspond to evaluating the weak form by restricting the test functions to the velocity and pressure components, respectively. Specifically, $(R_M)_{A,i}$ is the difference between the weak-form value $B(N_A\bm e_i,\,0;\,\bm u_h,p_h)$ and the right-hand side $F(N_A\bm e_i,\,0)$ when only the test function for the velocity is prescribed, whereas $(R_C)_A$ is the corresponding difference when only the test function for the pressure is prescribed. If $\bm R_M=\bm 0$ and $\bm R_C=\bm 0$, then the weak form is satisfied for all test functions, which means that the discretized equations are satisfied.

Consequently, the block matrices of the Jacobian $\bm J$ are expressed as
\begin{align}
  \bm J_{MM} &=
  \left.\frac{\partial \bm R_M}{\partial \bm u}\right|_{t+\Delta t,\,k},
  \label{eq:K_imp_clean}\\
  \bm J_{MC} &= \left.\frac{\partial \bm R_M}{\partial \bm p}\right|_{t+\Delta t,\,k},
  \label{eq:G_imp_clean}\\
  \bm J_{CM} &= \left.\frac{\partial \bm R_C}{\partial \bm u}\right|_{t+\Delta t,\,k},
  \label{eq:D_imp_clean}\\
  \bm J_{CC} &= \left.\frac{\partial \bm R_C}{\partial \bm p}\right|_{t+\Delta t,\,k}.
  \label{eq:L_imp_clean}
\end{align}

\newpage

%% file: include/appendixB.tex
\section{Hyper-Reduction Method}

\subsection{Hyper-Reduction via Empirical Cubature}
In G-ROMs, the computation of the projection operations for the Jacobian and residual vector, $\tilde{\bm J} = \bm{\Phi}^{\mathrm{T}}\bm{J}\bm{\Phi}$ and $\tilde{\bm r} = \bm{\Phi}^{\mathrm{T}}\bm{r}$, can often become a computational bottleneck. To reduce this computational cost, various hyper-reduction methods have been proposed, including DEIM (discrete empirical interpolation method) \cite{CHATURANTABUT20102737}, Gappy POD \cite{WILLCOX2006208,KANEKO2021104385}, GNAT (Gauss--Newton with approximated tensors) \cite{CARLBERG2013623}, S-OPT \cite{LAUZON2024B474}, ECSW (energy-conserving sampling and weighting) \cite{Farhat2015,GRIMBERG20211846}, and ECM (empirical cubature method) \cite{HERNANDEZ2017}. In this study, ECM \cite{HERNANDEZ2017} is employed as the hyper-reduction method, although the proposed method can be extended to other hyper-reduction methods.

Through finite element discretization, the analysis domain is divided into a set of elements $\mathcal{E}=\{e_1,\dots,e_{n_{\mathrm{elem}}}\}$. 
For each element $e$, a Boolean matrix $\bm{L}_e\in\mathbb{R}^{n_e\times n_{\mathrm{dof}}}$ is introduced to represent the correspondence between the local degrees of freedom defined within the element and the global degrees of freedom defined after element assembly. Here, $n_{\mathrm{dof}}$ denotes the total number of degrees of freedom in the FOM. If node $v_i$ has $d_i$ degrees of freedom, then
\[
  n_{\mathrm{dof}}=\sum_{i=1}^{n_{\mathrm{FOM}}}d_i .
\]
For systems with a constant number of degrees of freedom per node, this reduces to $n_{\mathrm{dof}}=d\,n_{\mathrm{FOM}}$, where $d$ is the number of degrees of freedom per node.
The Jacobian and residual vector can then be expressed in element-assembly form as
\begin{equation}
  \bm{r} = \sum_{e\in\mathcal{E}} \bm{L}_e^{\mathrm{T}} \bm{r}_e,
  \qquad
  \bm{J} = \sum_{e\in\mathcal{E}} \bm{L}_e^{\mathrm{T}} \bm{J}_e \bm{L}_e,
\end{equation}
where $\bm{r}_e\in\mathbb{R}^{n_e}$ is the residual vector associated with finite element $e$, $\bm{J}_e\in\mathbb{R}^{n_e\times n_e}$ is the Jacobian associated with finite element $e$, and $n_e$ is the number of degrees of freedom of finite element $e$.

In ECM, a subset $\tilde{\mathcal{E}}\subset\mathcal{E}$ consisting of $n_{\mathrm{HROM}}\ll n_{\mathrm{elem}}$ elements and nonnegative weights $w_e\ ({e\in\tilde{\mathcal{E}}})$ are determined, and the following approximations are used:
\begin{align}
  \bm{\Phi}^{\mathrm{T}} \bm{r}
  &= \sum_{e\in\mathcal{E}} \bm{\Phi}^{\mathrm{T}} \bm{L}_e^{\mathrm{T}} \bm{r}_e
   \;\approx\;
   \sum_{e\in\tilde{\mathcal{E}}} w_e\, \bm{\Phi}^{\mathrm{T}} \bm{L}_e^{\mathrm{T}} \bm{r}_e,
   \label{eq:ECM_r_B}\\
  \bm{\Phi}^{\mathrm{T}} \bm{J}\, \bm{\Phi}
  &= \sum_{e\in\mathcal{E}} \bm{\Phi}^{\mathrm{T}} \bm{L}_e^{\mathrm{T}} \bm{J}_e \bm{L}_e \bm{\Phi}
   \;\approx\;
   \sum_{e\in\tilde{\mathcal{E}}} w_e\, \bm{\Phi}^{\mathrm{T}} \bm{L}_e^{\mathrm{T}} \bm{J}_e \bm{L}_e \bm{\Phi},
   \label{eq:ECM_J_B}
\end{align}
where $\tilde{\mathcal{E}}$ and $w_e\ ({e\in\tilde{\mathcal{E}}})$ are determined by the singular value decomposition and the NNLS (nonnegative least-squares) problem described below.

In the training stage of ECM, $\tilde{\mathcal{E}}$ and $w_e$ are determined in a data-driven manner based on the solutions of the reduced system in \eqref{eq:NR_reduced} and \eqref{eq:NR_reduced_ansvec}. In the following, snapshot data $\{\bm{q}_s\}_{s=1}^{n_{\mathrm{snap}}}$, with $\bm{q}_s\in\mathbb{R}^{n_{\mathrm{POD}}}$, are prepared for the reduced problem \eqref{eq:NR_reduced_ansvec} in the POD-Galerkin method, and the high-dimensional solution vector in the ROM corresponding to the $s$-th snapshot data is taken as $\bm{\Phi}\bm{q}_s$.

First, for each element $e$, define
\begin{equation}
  \bm{g}_s^e
  := 
     \bm{\Phi}^{\mathrm{T}} \bm{L}_e^{\mathrm{T}} \bm{r}_e (\bm{L}_e\bm{\Phi} \bm{q}_s)
  \in\mathbb{R}^{n_{\mathrm{POD}}},
  \qquad
  \bm{b}_s := \sum_{e\in\mathcal{E}} \bm{g}_s^e \in\mathbb{R}^{n_{\mathrm{POD}}},
\end{equation}
where $\bm{g}_s^e$ is the vector associated with finite element $e$ and the $s$-th snapshot data, $s=1,\ldots,n_{\mathrm{snap}}$. Various datasets for storing $\bm{g}_s^e$ have been proposed, including those incorporating volume-constraint conditions \cite{YANO20192287}; in this study, the residual vector for each finite element in the above expression is stored as the dataset. For the unsteady diffusion equation in this study, $\bm{g}_s^e = \bm{\Phi}^{\mathrm{T}} \bm{L}_e^{\mathrm{T}} \bm{K}_e \bm{L}_e \bm{\Phi} \bm{q}_s - \bm{\Phi}^{\mathrm{T}} \bm{L}_e^{\mathrm{T}} \bm{f}_e$ was used. Here, $\bm{K}_e\in\mathbb{R}^{n_e\times n_e}$ is the element coefficient matrix associated with finite element $e$, and $\bm{f}_e\in\mathbb{R}^{n_e}$ is the element right-hand-side vector associated with finite element $e$. By stacking the above quantities vertically for all snapshots, the following matrix and vector are obtained:
\begin{equation}
\bm{G} = 
\begin{bmatrix}
\bm{g}_{1}^1 & \cdots & \bm{g}_{1}^{n_{\mathrm{elem}}} \\
\vdots  & \ddots & \vdots      \\
\bm{g}_{n_{\mathrm{snap}}}^1 & \cdots & \bm{g}_{n_{\mathrm{snap}}}^{n_{\mathrm{elem}}}
\end{bmatrix}
\in\mathbb{R}^{(n_{\mathrm{POD}}\,n_{\mathrm{snap}})\times n_{\mathrm{elem}}},
\qquad
\bm{b} = 
\begin{bmatrix}
\bm{b}_1 \\ \vdots \\ \bm{b}_{n_{\mathrm{snap}}}
\end{bmatrix}
\in\mathbb{R}^{n_{\mathrm{POD}}\,n_{\mathrm{snap}}}.
\end{equation}

In ECM, the singular value decomposition is applied to the matrix $\bm G$:
\begin{equation}
  \bm G
  =
  \bm Z
  \bm \Lambda
  \bm Y^{\mathrm T}.
\end{equation}
Based on the decay of the singular values, the dominant $p$ left singular vectors are extracted, yielding $\bm Z_p\in\mathbb{R}^{(n_{\mathrm{POD}}n_{\mathrm{snap}})\times p}$. This provides a compressed representation of $\bm G$ in the subspace spanned by a small number of orthogonal basis vectors.

Subsequently, the matrix $\bm G$ and vector $\bm b$ are projected onto the subspace compressed by the singular value decomposition, and the weight vector $\bm{w} = [w_1,\ldots,w_{n_{\mathrm{elem}}}]^{\mathrm{T}}$ is determined by
\begin{equation}\label{eq:ECM_global_opt}
\begin{gathered}
\bm{w}
=
\underset{\bm z\in\mathbb{R}^{n_{\mathrm{elem}}}}{\operatorname{arg\,min}}
\|\bm z\|_0
\\
\mathrm{subject\ to}
\quad
\left\|
\bm Z_p^{\mathrm T}
(\bm G\bm z-\bm b)
\right\|
\le
\varepsilon_{\mathrm{HROM}}\|\bm Z_{p}^{\mathrm T}\bm b\| ,
\quad
\bm z\succeq \bm 0 ,
\end{gathered}
\end{equation}
where $\|\cdot\|_0$ denotes the $\ell_0$ pseudo-norm, that is, the number of nonzero components of its argument. The parameter $\varepsilon_\mathrm{HROM}$ is a user-defined tolerance. The condition $\bm z\succeq \bm 0$ indicates that each component of the candidate weight vector $\bm z$ is nonnegative.
Using the set of elements corresponding to the nonzero components of the solution $\bm{w}$ of Eq.~\eqref{eq:ECM_global_opt},
\begin{equation}\label{ecm_sparse_elem}
  \tilde{\mathcal{E}} := \{\, e \in \mathcal{E} \mid w_e > 0 \,\}, \qquad
  n_{\mathrm{HROM}} := |\tilde{\mathcal{E}}| \ll n_{\mathrm{elem}},
\end{equation}
the online computation is accelerated using Eqs.~\eqref{eq:ECM_r_B} and \eqref{eq:ECM_J_B}.

For Eq.~\eqref{eq:ECM_global_opt}, finding the solution with optimal sparsity, that is, the minimum number of elements, is known to be NP-hard. Therefore, in this study, a greedy method is used to determine the sparse weights $\bm{w}$ and element set $\tilde{\mathcal{E}}$. Specifically, greedy NNLS \cite{LAWSON1995}, a greedy iterative algorithm, is adopted.

Because the proposed method does not depend on the choice of hyper-reduction method, it can also be extended to other methods such as DEIM, GNAT, S-OPT, and ECSW.

\subsection{Domain-Decomposition Hyper-Reduction}

\begin{figure*}[t]
    \begin{center}
    \hspace{25pt}
    \includegraphics[scale=0.65]{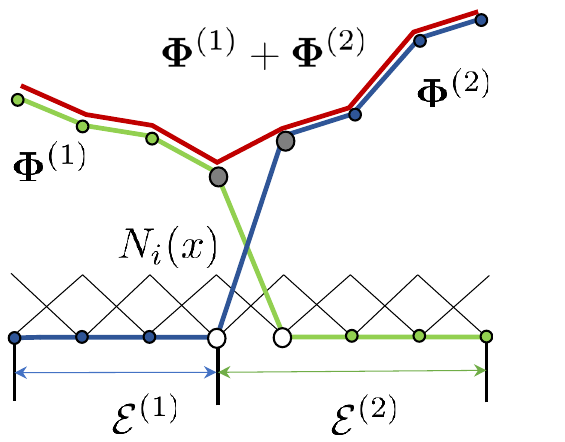}
    \end{center}
    \vspace{-10pt}
    \caption{\baselineskip=13pt Schematic overview of local POD bases and the subdomain-wise element partitioning used for ECM-based hyper-reduction.}
    \label{fig:LPOD_ECM}
\end{figure*}

This subsection presents a procedure for hyper-reduction based on the ECSW method of Farhat et al.~\cite{Farhat2015}, in which NNLS is applied locally to each subdomain to determine the weights. In this study, as in ECM, singular value decomposition is applied to the element-contribution matrix constructed in each subdomain, and a local NNLS problem is solved on the resulting low-dimensional basis. Because POD is also performed for each subdomain in this study, the POD basis matrix has a block-diagonal form as in Eq.~\eqref{eq:localPODmodes}.

Let the set of all elements be $\mathcal{E}=\{e_1,\ldots,e_{n_{\mathrm{elem}}}\}$. For each subdomain $i\;(1\le i\le N_{\mathrm{POD}})$, introduce an overlapping element set $\overline{\mathcal{E}}^{(i)}\subset\mathcal{E}$ based on the overlapping node set $\bar V^{(i)}$, such that
\begin{equation}
  \mathcal{E} \;=\; \bigcup_{i=1}^{N_{\mathrm{POD}}} \overline{\mathcal{E}}^{(i)}.
\end{equation}
Because elements may be shared among subdomains, if subdomain $i$ is adjacent to subdomain $j$, then
\begin{equation}\overline{\mathcal{E}}^{(i)} \cap \overline{\mathcal{E}}^{(j)} \neq \emptyset\label{eq:overlap_cover}\end{equation}
may hold.
To avoid double counting of elements during the hyper-reduction computation, a unique owner subdomain $\operatorname{own}(e)$ is assigned to each element $e$ as
\begin{equation}
  \operatorname{own}(e)
  := \max\bigl\{\, i \ \big| \ e\in\overline{\mathcal{E}}^{(i)} \,\bigr\}.
  \label{eq:owner}
\end{equation}
This rule is adopted for reproducibility in defining the algorithm.

Next, by defining the non-overlapping element set
\begin{equation}
  \mathcal{E}^{(i)}
  := \bigl\{\, e\in\mathcal{E} \ \big| \ \operatorname{own}(e)=i \,\bigr\},
  \label{eq:Ei_def}
\end{equation}
the collection $\{\mathcal{E}^{(i)}\}_{i=1}^{N_{\mathrm{POD}}}$ gives a mutually disjoint partition of the elements:
\begin{align}
  \mathcal{E}^{(i)} \cap \mathcal{E}^{(j)} &= \emptyset \quad (i\neq j), \label{eq:E_disjoint}\\
  \bigcup_{i=1}^{N_{\mathrm{POD}}}\mathcal{E}^{(i)} &= \mathcal{E}. \label{eq:E_cover}
\end{align}
Thus, the following relations are satisfied:
\begin{align}
  \bm{r}
  &= \sum_{i=1}^{N_{\mathrm{POD}}} \sum_{e\in\mathcal{E}^{(i)}}
     {\bm{L}_e}^{\mathrm{T}}\,\bm{r}_e,
  \label{eq:r_assembly}\\
  \bm{J}
  &= \sum_{i=1}^{N_{\mathrm{POD}}} \sum_{e\in\mathcal{E}^{(i)}}
     {\bm{L}_e}^{\mathrm{T}}\,\bm{J}_e\,\bm{L}_e.
  \label{eq:J_assembly}
\end{align}
Figure~\ref{fig:LPOD_ECM} shows the relationship between the L-POD basis and the element partition.
For each subdomain $i$, the elements in the set $\mathcal E^{(i)}$ are assigned local indices and expressed as
\begin{equation}
  \mathcal E^{(i)}
  =
  \left\{
  e^{(i)}_1,\ldots,e^{(i)}_{n_{\mathrm{elem}}^{(i)}}
  \right\},
  \label{eq:Ei_local_index}
\end{equation}
where $e^{(i)}_\ell\in\mathcal E^{(i)}$ denotes the $\ell$-th element belonging to subdomain $i$, and the index $\ell$ is the local element number within subdomain $i$.

In domain-decomposition ECM, the global matrix $\bm{G}$ is partitioned as
\begin{equation}
  \bm{G} = \bigl(\bm{G}^{(1)},\ldots,\bm{G}^{(N_{\mathrm{POD}})}\bigr),
\end{equation}
where the matrix $\bm G^{(i)}$ and vector $\bm b^{(i)}$ are defined as
\begin{align}
\bm{G}^{(i)} &=
\begin{bmatrix}
\bm g_{1}^{e^{(i)}_1,(i)} & \cdots & \bm g_{1}^{e^{(i)}_{n_{\mathrm{elem}}^{(i)}},(i)}\\
\vdots & \ddots & \vdots\\
\bm g_{n_{\mathrm{snap}}}^{e^{(i)}_1,(i)} & \cdots & \bm g_{n_{\mathrm{snap}}}^{e^{(i)}_{n_{\mathrm{elem}}^{(i)}},(i)}
\end{bmatrix}
\in
\mathbb{R}^{\left(n_{\mathrm{snap}} \sum_{j=1}^{N_{\mathrm{POD}}} n_{\mathrm{POD}}^{(j)}\right)\times n_{\mathrm{elem}}^{(i)}},
\\
\bm b^{(i)} &=
\begin{bmatrix}
\bm b_{1}^{(i)}\\
\vdots\\
\bm b_{n_{\mathrm{snap}}}^{(i)}
\end{bmatrix}
\in
\mathbb{R}^{\,n_{\mathrm{snap}} \sum_{j=1}^{N_{\mathrm{POD}}} n_{\mathrm{POD}}^{(j)}}.
\label{eq:Gibidef}
\end{align}
For an element $e\in\mathcal{E}^{(i)}$ in subdomain $i$, the components are defined as
\begin{equation}
  \bm g_{s}^{e,(i)}
  :=
  \bm \Phi^{\mathrm{T}} \bm L_e^{\mathrm{T}} \bm r_{e} (\bm L_{e} \bm{\Phi} \bm q_s)
  \;\in\; \mathbb{R}^{\sum_{j=1}^{N_{\mathrm{POD}}}n_{\mathrm{POD}}^{(j)}},
  \qquad
  \bm b_s^{(i)} := \sum_{e\in\mathcal{E}^{(i)}} \bm g_{s}^{e,(i)} \in \mathbb{R}^{\sum_{j=1}^{N_{\mathrm{POD}}}n_{\mathrm{POD}}^{(j)}}.
  \label{eq:local_gb}
\end{equation}
Here, the index corresponding to the snapshot data is $s=1,\ldots,n_{\mathrm{snap}}$, and the reduced solution vector is $\bm{q}_s\in\mathbb{R}^{\sum_{i=1}^{N_{\mathrm{POD}}} n_{\mathrm{POD}}^{(i)}}$. In practical computations, $\bm G^{(i)}$ and $\bm b^{(i)}$ are computed using the POD bases of subdomain $i$ and its neighboring subdomains ${j}\, \in\, \mathrm{nbhd}(\Tilde{v}_i, \Tilde{G})$, and therefore the computation can be performed independently for each subdomain.

Next, singular value decomposition is applied to the matrix $\bm G^{(i)}$ constructed in each subdomain $i$:
\begin{equation}
  \bm G^{(i)}
  =
  \bm Z^{(i)} \bm \Lambda^{(i)}
  \left(\bm Y^{(i)}\right)^{\mathrm T}.
  \label{eq:local_svd}
\end{equation}
Based on the decay of the singular values, the dominant $p_i$ left singular vectors $\bm Z_{p_i}^{(i)}$ are extracted:
\begin{equation}
  \bm Z_{p_i}^{(i)}
  =
  \left[
  \bm z_1^{(i)},\ldots,\bm z_{p_i}^{(i)}
  \right]
  \in
  \mathbb{R}^{\left(n_{\mathrm{snap}} \sum_{j=1}^{N_{\mathrm{POD}}} n_{\mathrm{POD}}^{(j)}\right)\times p_i}.
  \label{eq:local_svd_basis}
\end{equation}
This represents $\bm{G}^{(i)}$ in each subdomain using a small number of orthogonal basis vectors.

Finally, in each subdomain $i$, the following optimization problem in the space compressed by singular value decomposition is solved using greedy NNLS \cite{LAWSON1995}:
\begin{equation}\label{eq:ECM_local_opt}
\begin{gathered}
\bm w^{(i)}
=
\underset{\bm z^{(i)}\in\mathbb R^{n_{\mathrm{elem}}^{(i)}}}
{\operatorname{arg\,min}}
\|\bm z^{(i)}\|_0
\\
\mathrm{subject\ to}
\quad
\left\|
\left(\bm Z_{p_i}^{(i)}\right)^{\mathrm T}
(\bm G^{(i)}\bm z^{(i)}-\bm b^{(i)})
\right\|
\le
\varepsilon_{\mathrm{HROM}}\|\left(\bm Z_{p_i}^{(i)}\right)^{\mathrm T}\bm b^{(i)}\| ,
\quad
\bm z^{(i)}\succeq 0.
\end{gathered}
\end{equation}
The resulting solution is the local weight vector $\bm w^{(i)}\in\mathbb{R}^{n_{\mathrm{elem}}^{(i)}}$.